\documentclass[titlepage,12pt]{amsart}
\usepackage[stable]{footmisc}
\makeatletter
\def\myfnt{\ifx\protect\@typeset@protect\expandafter\footnote\else\expandafter\@gobble\fi}
\makeatother

\usepackage{amsfonts}
 \usepackage{amssymb}
 \usepackage{amsmath,amsxtra,amsthm}

\usepackage[usenames]{color}

\usepackage[dvips]{graphicx}

\usepackage[all,matrix,arrow,curve]{xy}

\usepackage{mathrsfs}

\DeclareFontFamily{OT1}{pzc}{}
\DeclareFontShape{OT1}{pzc}{m}{it}{<-> s * [1.15] pzcmi7t}{}
\DeclareMathAlphabet{\mathpzc}{OT1}{pzc}{m}{it}

\DeclareSymbolFont{bbold}{U}{bbold}{m}{n}
\DeclareSymbolFontAlphabet{\mathbbold}{bbold}

\newtheorem{theorem}{Theorem}
\newtheorem{deff}{Definition}

\newtheorem{proposition}{Proposition}
\newtheorem{example}{Example}
\newtheorem{lemma}{Lemma}

\newtheorem{cor}{Corollary}

\newtheorem{rem}{Remark}

\renewcommand{\proof}{{\bf Proof.}~}

\newcommand{\mto}{\mapsto}
\newcommand{\bqa}{\begin{eqnarray}}
\newcommand\eqa {\end{eqnarray}}
\newcommand{\beq}{\begin{eqnarray}}
\newcommand{\beqn}{\begin{eqnarray}\nonumber}
\newcommand{\eeq}{\end{eqnarray}}
\newcommand{\be}{\begin{array}}
\newcommand{\ee}{\end{array}}

 \newcommand{\pt}{\partial}

   \newcommand\vf\varphi

 \newcommand{\Hom}{\mathrm{Hom}}
 \newcommand{\End}{\mathrm{End}}
 \newcommand{\Ker}{\mathrm{Ker}}
 \newcommand{\Id}{\mathrm{Id}}

 \newcommand{\rk}{\mathrm{rk}}
 \newcommand{\md}{\mathrm{d}}

 \newcommand{\D}{\mathrm{D}}

 \newcommand{\cP}{{\mathcal{P}}}
 \newcommand{\cC}{{\mathcal C}}
 \newcommand{\cAfeld}{{\mathcal{A}}}

 \newcommand{\cL}{\mathcal{L}}
 
 \newcommand{\cG}{{\mathcal{G}}}
 
 \newcommand{\cF}{{\mathcal{F}}}

 \newcommand{\R}{{\mathbb R}}

 \newcommand{\N}{{\mathbb N}}

  \def\g{{\mathfrak g}}

  \def\h{{\mathfrak h}}

   \def\a{\alpha}

   \def\de{\delta}

  \def\sst{\scriptscriptstyle}

\title[Lie algebroids with compatible geometrical structures]{Lie algebroids, gauge theories, \\and  compatible geometrical structures}%
\author{Alexei Kotov}
\address{Alexei Kotov: Faculty of Science, University of Hradec Kralove, Rokitanskeho 62, Hradec Kralove
50003, Czech Republic}
\email{oleksii.kotovATuhk.cz}

\author{Thomas Strobl}
\address{Thomas Strobl: Institut Camille Jordan, Universit\'e Claude Bernard Lyon 1, Universit\'e de Lyon,
43 boulevard du 11 novembre 1918, 69622 Villeurbanne Cedex, France}
\email{stroblATmath.univ-lyon1.fr}

\date{March 7, 2018}
\begin{document}

\begin{abstract}
The construction of gauge theories beyond the realm of Lie groups and algebras leads one to consider Lie groupoids and algebroids equipped with additional geometrical structures which, for gauge invariance of the construction, need to satisfy particular compatibility conditions. This paper is supposed to analyze these compatibilities from a mathematical perspective.

\vskip 2mm\noindent
In particular, we show that the compatibility of a connection with a Lie algebroid that one finds is the Cartan condition, introduced previously by A.\ Blaom. For the metric on the base $M$ of a Lie algebroid equipped with any connection, we show that the compatibility suggested from gauge theories implies that the foliation induced by the Lie algebroid becomes a Riemannian foliation. Building upon a result of del Hoyo and Fernandes, we prove furthermore that every Lie algebroid integrating to a proper Lie groupoid admits a compatible Riemannian base. We also consider the case where the base is equipped with a compatible symplectic or  generalized
Riemannian structure.

\vspace{20mm}
\noindent \keywordsname : Lie algebroids, Riemannian foliations, symplectic realization,
generalized geometry, symmetries and reduction, gauge theories.

\vspace{2.5cm}

\tableofcontents
\end{abstract}
\setcounter{page}{0}

\maketitle


\clearpage

\setcounter{page}{1}

\section{Introduction}

\vskip 3mm\noindent Lie algebroids were introduced in the 1960s by Pradines \cite{Pradines67} as a formalization of ideas going back to works of Lie and Cartan. By definition it is a Lie algebra structure on the sections of a vector bundle $A\to M$ which in turn is equipped with a bundle map $\rho \colon A \to TM$ such that for every pair of sections $s,s' \in \Gamma(A)$ and every function $f\in C^\infty(M)$ one has the Leibniz rule $[s,fs']=f[s,s']+(\rho(s)f) s'$. They combine usual geometry and Lie algebra theory below a common roof, interpolating between tangent bundles and foliation distributions on the one hand and Lie algebras and their actions on manifolds on the other hand. Correspondingly, Lie algebroids provide a particular way to generalize notions of ordinary geometry to geometry on them: Typical examples include the generalization of the de Rham differential, ${}^A\md \colon \Gamma(\Lambda^\bullet A^*) \to \Gamma(\Lambda^{\bullet + 1} A^*)$, Lie algebroid covariant derivatives ${}^A{}\nabla$, where sections $\psi$ of another vector bundle $V$ over the same base can be differentiated along sections $s$ of $A$, ${}^A{}\nabla_s \psi \in \Gamma(V)$, together with a notion of its curvature or torsion, or the fact that any fiber metric ${}^A g$ on $A$ induces a unique $A$-torsion-free, metrical $A$-connection on $A$. Likewise, notions from standard Lie theory generalize to this setting:
For example, representations of a Lie algebra on a vector space are replaced here by flat $A$-connections on vector bundles.\footnote{We will assume acquaintance to Lie algebroids and Lie groupoids throughout this paper, see, e.g., \cite{daSilva-Weinstein99,Mackenzie05,Moerdijk-Mrcun03} for the necessary background material.}


\vskip 3mm\noindent The purpose of the present paper is to investigate the interplay between Lie algebroids and additional (ordinary) geometrical structures defined on them. These structures together with the appropriate compatibilities are inspired by mathematical physics, where they appear naturally in the context of gauge theories, see~\cite{Kotov-Strobl14,Kotov-Strobl_CYMH15} for the corresponding papers in physics, or the more recent papers \cite{withyoungpeople,Constrained}, written in a language for a mixed audience, as well as the subsequent section, {\bf Sec.~\ref{Motivation}}, for a summary in a purely mathematical language. There are different options of interest in this context. All of them have in common that $A$ is equipped with an ordinary connection $\nabla$.

 \vskip 3mm\noindent One of the natural questions posing itself is thus the one of a good compatibility of $\nabla$ with a Lie algebroid structure on $A$.
An evident option is that $\nabla$ should preserve the bracket. This is too restrictive, however:
 Suppose, for example, that $A=T^*M$ for a Poisson manifold $(M,\Pi)$.\footnote{For simplicity, we sometimes write simply $A$ for the Lie algebroid data $(A,\rho,[\cdot,\cdot])$. Depending on the context, $A$ either denotes a Lie algebroid or only its underlying vector bundle.} Then such a condition corresponds to $\nabla \Pi =0$,  restricting $\Pi$ to have constant rank and thus excluding most of the interesting examples of Poisson manifolds. Another natural option is the following one: Any ordinary connection on $A$ gives rise to a very trivial $A$-connection by ${}^A{}\nabla_\cdot = \nabla_{\rho(\cdot)}$. This does not yet contain information about the Lie bracket on $A$. However, there is a notion of a duality of $A$-connections on a
Lie algebroid (cf.~\cite{Blaom04}),
${}^A{}\nabla \mapsto {}^A{}\nabla^*$, where ${}^A{}\nabla^*_s s'= [s,s']+{}^A{}\nabla_{s'}s$. One then may ask the $A$-curvature of ${}^A{}\nabla^*$ to vanish.
This is too weak now, however, since for a bundle of Lie algebras, $A$ with $\rho \equiv 0$, this condition becomes vacuous, where here we would very much like to preserve the bracket.

\vskip 3mm\noindent The good notion is the one that was introduced by A.\ Blaom in  \cite{Blaom04,Blaom05}: Any connection in $A$ gives a splitting of $J^1(A) \to A$, where $J^1(A)$ denotes the first jet bundle of sections of $A$. Since $J^1(A)$ carries a natural Lie algebroid structure induced by $A$, we require the $\nabla$-induced splitting to be a Lie algebroid morphism. This condition implies indeed that ${}^A{}R_{{}^A{}\nabla^*}=0$ and still reduces to the wished-for constancy of the bracket for bundles of Lie algebras. We will recall this notion in detail and show in particular that it reproduces the mathematically less transparent formulation that has been found in the context of mathematical physics. We call a compatible couple $(A,\nabla)$ a Cartan-Lie algebroid.

\vskip 3mm\noindent In the applications, the bundle $A$ will then be equipped with additional geometric structures, like a Riemannian metric $g$ on the base or a fiber metric ${}^A{}g$ on the bundle. In the present paper, we will focus merely on geometrical structures defined on $M$, starting in {\bf Sec.~\ref{section:Killing-Lie}} with Riemannian ones.
The compatibility that a metric $g$ on the base of a Lie algebroid $A$ is required to satisfy will be of the form \begin{equation} \label{one} {}^\tau{}\nabla g = 0\,,\end{equation} where ${}^\tau{}\nabla$ is an $A$-connection on $TM$, intertwined by $\rho$ with ${}^A{}\nabla^*$ introduced above. Since the transition to $TM$ implies also replacing the Lie bracket on $A$ by the canonical bracket of vector fields, the condition \eqref{one} does in fact only depend on the anchor $\rho \colon A \to TM$ and thus can be formulated for any anchored bundle already. While some of the statements in this paper will be formulated for this enlarged setting, one of the main results of the present paper requires a Lie algebroid that is moreover integrable to a proper Lie groupoid. Properness here means that the map $s\times t \colon \cG \times \cG \to M$ is proper, i.e.~compact subsets of $M$ have compact pre-images. Thus these are Lie algebroids which  generalize compact Lie algebras. Building on a result of del Hoyo and Fernandes \cite{delHoyo-Fernandes15}, we show that any such a Lie algebroid admits a metric $g$ satisfying Eq.~\eqref{one}.

\vskip 3mm\noindent We call a metric $g$ compatible with a Lie algebroid $A$ in the sense of Eq.~\eqref{one} a Killing Lie algebroid. Although this name is difficult to generalize to other geometrical structures, it expresses the close relation of it to symmetries of the geometric structure on the base, here the metric: Indeed, the image with respect to the anchor of any (locally) covariantly constant section is rather easily seen to be a (local) Killing vector field of $g$. Be there Killing vectors of the metric $g$ on $M$ or not, the notion also guarantees good quotients whenever the leaf-space of the Lie algebroid permits a good quotient; in this case the natural projection map becomes a Riemannian submersion. A natural generalization of Riemannian submersions to singular foliations are Riemannian foliations \cite{Molino88}. We will show that the natural foliation of $M$ induced by any Killing Lie algebroid is a Riemannian foliation. If the reverse is true, i.e.~if any Lie algebroid $A$ equipped with a Riemannian metric $g$ on its base can be equipped with a connection $\nabla$ so that the data become those of a Killing Lie algebroid, is an interesting problem left open in this article. Some simple examples and facts about Killing Lie algebroids where the connection $\nabla$ on $A$ is flat, are deferred to an appendix, {\bf App.~A}.

\vskip 3mm\noindent In {\bf Sec.~\ref{section:symplectic}} we briefly consider other geometrical structures on the base, without developing this subject in the present paper to much of a depth. The reason for this is two-fold: First, while also for other compatible geometrical structures on the base of  a Lie algebroid with a smooth leaf space $Q=M/\cF$ permit the geometric structure to descend to the quotient in some way, the resulting structure on $Q$ can be of a different type. We will illustrate this with a symplectic structure, where on the quotient one induces a Poisson structure, so that $(M,\omega) \to (Q,\Pi_Q)$ provides a symplectic realization of this quotient. The second reason is that some considerations in string theory \cite{withyoungpeople} suggest a more general compatibility condition than Eq.~\eqref{one} induced by a single connection $\nabla$ on $A$. We want to show how these structures can be hosted as well in the present framework and related in this particular context to what is called now generalized geometry \cite{Hitchin03}.

\vskip 3mm\noindent There are two related studies that are treated in accompanying papers of us, which are worth mentioning here. First, as mentioned above, the notion of a compatibility of a metric $g$ on the base $M$ of a vector bundle $A$ is already well-defined, when $A$ is an anchored bundle only, i.e.\ a bundle $A$ equipped with a map $\rho \colon A \to TM$ covering the identity map without any further conditions. In \cite{biinvariant} we show that one can extend this in a unique way to an infinite-rank free Cartan Killing Lie algebroid $(FR(A),\widetilde{\nabla},g)$. Here $FR(A)$ is the free Lie algebroid of the anchored bundle $A$ as introduced by A.\ Kapranov \cite{Kapranov07}, $\widetilde{\nabla}$ the unique extension of $\nabla$ to a Cartan-compatible connection on the Lie algebroid $FR(E)$, and $g$ the original metric on $M$, which now is also compatible with this enlagred structure.

\vskip 3mm\noindent Second, if one has a Cartan-Lie algebroid $(A,\nabla)$ together with a metric $g$ on its base, compatible in the sense of equation \eqref{one}, and a fiber metric ${}^Ag$, compatible with $(A,\nabla)$ in the sense of ${}^A{}\nabla^* {}^Ag=0$ (with the $A$-connection introduced above), one may see that these data give the appropriate generalisation of a quadratic Lie algebra to the setting of Lie algebroids. In \cite{universal} we determine the necessary and sufficient conditions for this to integrate to a Riemannian groupoid for the case that the Lie algebroid $A$ itself is integrable to a groupoid. Interestingly, there are obstructions to this integration in general.

\vskip 5mm \noindent {\bf Acknowledgements.} We are grateful to  Matias del Hoyo, Rui Fernandes,
Camille Laurent-Gengoux, Llohann Speran\c{c}a and  Alan Weinstein for their interest and useful discussions. Likewise we thank Anton Alekseev, Jim Stasheff and Marco Zambon for remarks on the manuscript. 
The research of A.K. was supported by grant no. 18-00496S of the Czech Science Foundation. The research of T.S. was supported the project MODFLAT of the European Research Council (ERC) and the NCCR SwissMAP of the Swiss National Science Foundation. We both gratefully acknowledge partial support by Projeto P.V.E.~88881.030367/2013-01 (CAPES/Brazil).

\section{Motivation via Mathematical Physics}

\label{Motivation}

\vskip 3mm\noindent
Physics often points to new, mathematically interesting notions, worth being studied in their own right, and selecting the ``right choice'' out of a variety of otherwise conceivable options. The present article implicitly pursues such a study, suggested by a generalization of gauge theories of the Yang-Mills type where the structure  group or Lie algebra is replaced by a structural Lie groupoid/algebroid. While keeping the main text below as well as eventual sequels to this article mathematically self-contained without it, we provide a short summary of this motivation here (for more details cf.~our previous, more physics oriented articles \cite{Kotov-Strobl14,Kotov-Strobl_CYMH15} as well as \cite{Mayer-Strobl09,withyoungpeople}).

\vskip 3mm\noindent
Standard sigma models are defined after selecting a source (pseudo) Riemannian manifold $(\Sigma,h)$ and a target Riemannian manifold $(M,g)$. This gives rise to a functional $S_0[X]$ on the space of maps $X$ from $\Sigma$ to $M$, the Euler-Lagrange equations of which are harmonic maps. The isometry/symmetry group $G$ of  $(M,g)$
lifts to $S_0$. ``Gauging'' of this symmetry group $G$ then amounts to the following procedure:
One considers a $G$-principal bundle $P$ over $\Sigma$. For simplicity, we restrict to a trivial and trivialized bundle $P := \Sigma \times G$. Gauging is effectuated by introducing a connection in $P$ as an additional ``field'', which we may identify with a $\g$-valued 1-form ${\cAfeld}$
on $\Sigma$, where $\g = \mathrm{Lie}(G)$. Replacing the tangent map $\mathrm{d} X$ of $X$ by its covariantization
$\mathrm{D}X := \mathrm{d} X - \rho({\cAfeld})$ in the original functional $S_0$, where the bundle map $\rho \colon M\times \g \to TM$ encodes the $\g$-action on $M$, one obtains a new functional $S_1[X,{\cAfeld}]$. This ``gauged theory'' $S_1$ has the enhanced symmetry group $C^\infty(\Sigma,G)$, which is called the group of gauge transformations.

\vskip 3mm\noindent
Both of the ``fields'' together, the maps $X \in C^\infty(\Sigma,M)$ called scalar fields, as well as the connection 1-forms ${\cAfeld} \in \Omega^1(\Sigma,\g)$, called gauge fields, can be reinterpreted as a vector bundle morphism  $a \colon T\Sigma \to A$, where $A=M \times \g$ may be viewed as an action Lie algebroid. The Killing equation expressing the isometry of $g$ on $M$ with respect to the Lie algebra $\g$ can be rewritten as follows: Consider the action Lie algebroid $(A, \rho, [\cdot, \cdot])$ over the Riemannian manifold $(M,g)$ together with the canonical flat connection $\nabla$ on $A$ such that
\begin{equation}\label{barrho}
\left(\nabla (\iota_\rho g )\right)_{sym} = 0.
\end{equation}
Here $\iota_\rho g = g(\rho, \cdot) \in \Gamma(A^* \otimes T^*M)$ and for the $TM$-part of $\nabla$ one uses the Levi-Civita connection of $g$. The main observation coming from mathematical physics of relevance for the present article is that one can introduce a gauged theory of scalar fields for the triple $(A,\nabla,g)$ of a \emph{general} Lie algebroid $A$ equipped with a connection $\nabla$ and defined over a Riemannian base $(M,g)$ \emph{provided} the compatibility equation \eqref{barrho} is satisfied \cite{Kotov-Strobl14}. We call such a triple a Killing Lie algebroid, since for any constant section $s \in \Gamma(A)$, $\nabla s = 0$, the vector field $\rho(s)$ generates an infinitesimal isometry of $g$.

\vskip 3mm\noindent
Let us be more specific here. ``Gauging'' in this generalized context, i.e., in particular, without necessarily isometries of $g$, applies to the following situation: One is given a functional $S_0$ depending on $X \colon \Sigma \to M$ as described above. Assume that there is a Lie algebroid $(A,\rho,[ \cdot , \cdot])$ over $M$ equipped with a connection such that Eq.~\eqref{barrho} is satisfied. In fact, for what follows it is sufficient to have an anchored bundle $(A,\rho)$ over $M$, with an integrable singular distribution $\rho(A) \subset TM$. It is only important that Eq.~\eqref{barrho} holds true.
Then there exists a functional $S_1$ depending on $a \colon T\Sigma \to A$ with the following properties:
\begin{enumerate}
\item $S_1$ restricted to such maps $a$ that map to the 0-section in $A$ coincides with $S_0$.
\item Vector fields along the foliation of $M$ tensored by functions over $\Sigma$ can be lifted to infinitesimal symmetries of $S_1$ (gauge transformations).
\item
Whenever the quotient space $Q=M/\!\sim$ is smooth and $\Sigma$ is contractible, the equivalence classes of solutions to the Euler-Lagrange equations of $S_1$ modulo gauge symmetries are in one-to-one correspondence with solutions of an ``ungauged'' theory $S_0$ where $(M,g)$ is  replaced by the quotient $(Q,g_Q)$.
\end{enumerate}
In the last item, the leaf space is given by the  orbits generated by $\rho(A)$; Eq.~\eqref{barrho} then ensures that for smooth quotient $Q$ there is a unique metric $g_Q$ such that the projection from $M$ to $Q$ is a Riemannian submersion \cite{Kotov-Strobl14}. We also remark that the infinite-dimensional group of gauge symmetries is essentially defined already once the functional $S_1$ is given,   cf., e.g., \cite{Henneaux-Teitelboim}. While the third condition or property makes sense only for smooth quotients $M/\!\sim$, the gauge theory is particularly interesting or useful precisely when this quotient is not smooth (otherwise one would not need the concept of the gauge theory).

\vskip 3mm\noindent
Let $d$ be the dimension of the source manifold $\Sigma$, then the functional $S_0$ can be extended naturally, if one is given some $B \in \Omega^d(M)$,
namely, by adding its pullback by $X \colon \Sigma \to M$ to the original functional. In this way, by gauging one is led to consider $(A,\nabla,g,B)$, which, in addition to Eq.~\eqref{barrho}, has to satisfy
\begin{equation} \label{LieB}
\left(\iota_\rho \circ \mathrm{D} + \mathrm{D} \circ \iota_\rho \right) B = 0 \, .
\end{equation}
Here $\iota_\rho$ denotes contraction with the $TM$-part of $\rho \in \Gamma(A^* \otimes TM)$ and $\mathrm{D}$ is the exterior covariant derivative induced by $\nabla$. The similarity of the two equations \eqref{barrho} and \eqref{LieB} becomes more transparent, when rewritten in the following way:
\begin{eqnarray}
{\cL}_{\rho(s)} g &=& 2 \,\mathrm{Sym} \left( \iota_{(\mathrm{id} \otimes \rho) \nabla(s)} g \right), \label{Lie1}\\ {\cL}_{\rho(s)} B &=& d \, \mathrm{Alt} \left( \iota_{(\mathrm{id} \otimes \rho)\nabla(s)} B \right),
\label{Lie2}
\end{eqnarray}
which has to hold true for arbitrary sections $s\in \Gamma(A)$. Here on the r.h.s.~the map $\rho \colon \Gamma(A) \to \Gamma(TM)$ is applied to the second factor in $\nabla(s) \in \Gamma(T^*M \otimes A)$ and $\mathrm{Sym}$ and $\mathrm{Alt}$ denote the symmetrization and antisymmetrization projectors of the tensor product, respectively.
Again, for constant sections $s \in \Gamma(A)$, the vector fields $\rho(s)$ generate symmetries of $B$.

\vskip 3mm\noindent
Eq.~\eqref{LieB} can be used to motivate definitions in the context of (higher pre-) symplectic structures and reductions, and in fact so in potentially two different ways. We can either identify the (higher) symplectic form $\omega$ directly with $B$. Then,
Eq.~\eqref{LieB} suggests the form of an ``covariantized (higher) symplectic equation'' as $[\mathrm{D}, \iota_\rho] \omega = 0$, or, $\omega$ being closed, as
\begin{equation}\label{symp}
\D(\iota_\rho \omega) = 0 \, .
\end{equation}
Note that for $\nabla$ flat and $e_a$ a covariantly constant basis, this equation reduces to ${\cL}_{\rho_a} \omega = 0$, stating that the collection of vector fields $\rho_a$ is symplectic and higher symplectic for $d=2$ and $d>2$, respectively. This is e.g.~the case for the action Lie algebroid with its canonical flat connection. One certainly may also drop the non-degeneracy condition on $\omega$, replacing ``symplectic'' by ``pre-symplectic'' in the above.

\vskip 3mm\noindent Another option, more natural from the sigma model perspective, is to identify $\omega$ with $H=\mathrm{d}B$. Eq.~\eqref{LieB} then is best reinterpreted as suggesting a ``covariantized (higher) Hamiltonian equation'' of the form:
\begin{equation}\label{Ham}
\iota_\rho \omega = \D h \, ,
\end{equation}
where $\omega$ is a (higher) symplectic form,
$\omega \in \Omega_{cl}^{d+1}(M)$ non-degenerate, and $h \in \Omega_{cl}^{d-1}(M,A)$ a ``covariantized (higher) Hamiltonian''. For $d=1$, $(M,\omega)$ is a symplectic manifold in the usual sense  and the section $h \colon M \to A^*$ a ``covariantized moment map''.
Again, this reduces to the ordinary definition of a moment map $h \colon M \to \g^*$ for the case of $A=M\times \g$. More generally, whatever $A$ and $\nabla$, if one has a constant section $s$, then $\rho(s)$ is a Hamiltonian vector field with Hamiltonian $h(s)\equiv \langle h,s\rangle$. Note that, in contrast to standard symplectic geometry, the covariantized Hamiltonian equation, Eq.~\eqref{Ham}, does not imply and is not a special case of the covariantized symplectic equation, Eq.~\eqref{symp} (except for zero curvature $\mathrm{D}^2=0)$. In the present article we will consider the covariantized symplectic equation, postponing the Hamiltonian version to possible later work.\footnote{A.~Weinstein informed us, on the other hand, that he and C.~Blohmann have work in progress on the Hamiltonian case with $d=1$. They call $(A,\nabla,\omega)$ for $d=1$ satisfying Eq.~\eqref{Ham} \emph{Hamiltonian Lie algebroids}.}


\vskip 3mm\noindent There is yet another option of interpreting the \emph{couple} of equations \eqref{Lie1} and \eqref{Lie2} for $d=2$ of relevance for the  present paper. This is related to the notion of a generalized metric, cf., e.g., \cite{Hitchin03}: An ordinary metric $g$ on $M$ together with a not further restricted 2-form $B\in \Omega^2(M)$ can be viewed as a generalized metric $\Phi = g + B$ on $M$. It was seen in several instances already that strings are intimately related to generalized geometries, cf, e.g., \cite{Zabzineetal, Zabzine, Alekseev-Strobl, Kotov-Salnikov-Strobl14}. This corresponds to a theory with a worldsheet $\Sigma$ of two dimensions, $d=2$. It is also precisely this dimension where the Hodge dual of a 1-form gauge field ${\cAfeld}$ is again a 1-form. Thus, in $d=2$, and only there, the equations \eqref{Lie1} and \eqref{Lie2} admit a generalization, corresponding precisely to the additional freedom in the choice for gauge transformations (cf.~\cite{withyoungpeople} for the details). While given a Lie algebroid $A$ over $M$, the compatibility for $g$ and $B$ required in Eq.~\eqref{Lie1} and Eq.~\eqref{Lie2} is parameterized by a connection $\nabla$ in $A$, now there is in addition  a section $\psi \in \Gamma(T^*M \otimes \End(A))$ in the correspondingly weakened conditions, which then read as follows:
\begin{eqnarray}
{\cL}_{\rho(s)} g &=& 2\, \mathrm{Sym} \left( \iota_{(\mathrm{id} \otimes \rho) \nabla(s)} g +\iota_{(\mathrm{id} \otimes \rho)\psi(s)} B \right), \label{Lie3}\\ {\cL}_{\rho(s)} B &=& 2 \,\mathrm{Alt} \;\, \left( \iota_{(\mathrm{id} \otimes \rho)\nabla(s)} B \pm \iota_{(\mathrm{id} \otimes \rho) \psi(s)} g\right),
\nonumber
\end{eqnarray}
to hold true for all sections $s \in \Gamma(A)$. The sign in the last term in the second line corresponds to the signature of the metric $h$ on $\Sigma$: one has a plus sign, if $h$ has a Lorentzian signature and a minus sign, if the string theory is of Euclidean signature. In the main text below, we will interpret the coupled equations \eqref{Lie3} geometrically, using ideas from the realm of generalized geometry on $M$.

\vskip 3mm\noindent
Except for maybe this last, more involved system  \eqref{Lie3}, the guiding principle behind the notion one is tempted to introduce by the study of gauge theories should be clear by now: One always has a Lie algebroid $A$ with a connection $\nabla$. The base $M$ of $A$ is equipped with a geometrical structure ${\cG}=\{ g, \omega, \ldots \}$. Then there is a compatibility condition which implies that \emph{constant} sections generate symmetries of ${\cG}$ by means of the anchor $\rho$ (Eqs.~\eqref{barrho}, \eqref{symp} or \eqref{LieB}, \eqref{Ham}, etc). Moreover, in the case that $A=M\times \g$ is the action Lie algebroid with its canonical flat connection, one gets back the usual notions of Killing, (higher) symplectic, Hamiltonian, $\ldots$ vector fields or the notion of a moment map. But even if one does not deal with such symmetries of ${\cG}$ in the strict sense, the gauge theory formulation suggests that in such cases still a \emph{reduction} is meaningful---mathematically certainly with the additional assumption of a smooth quotient, or, more contemporary, interpreting the data as describing a quotient stack. i.e.~a particular smooth description of an otherwise singular quotient.

\vskip 3mm\noindent
Let us return to the logic immanent to gauge theories again. Physically one is compelled to add also a ``kinetic term'' for the gauge fields, at least if they should correspond to propagating interaction particles. In the standard setting with the isometry Lie algebra $\g$, one needs an ad-invariant metric on $\g$ in this context. In the generalization of $M \times \g$ to arbitrary Killing Lie algebroids $A$, this becomes a fiber metric ${}^A{}g$ on the vector bundle $A$ satisfying
\begin{equation} \label{Ag}
{}^A{}\widetilde{\nabla} \left({}^A{}g\right)=0
\end{equation}
with respect to the Lie algebroid or $A$-covariant derivative\footnote{In the body of the paper, this $A$-connection is mostly denoted by ${}^\alpha{}\nabla$, or by ${}^A{}\nabla^*$, so as to express that it is dual to the evident $A$-connection ${}^A{}\nabla_\cdot = \nabla_{\rho(\cdot)}$ in a precise way, cf.~Eq.\eqref{A-conn-dual} below. The $A$-torsion ${}^A{}\widetilde{T}$ of ${}^A{}\widetilde{\nabla}$ used below is the negative of the one of ${}^A{}\nabla$, denoted by ${}^A{}T$ in the main text; this difference in sign is important for comparison of formulas in this section with those in the rest of the paper.}
\begin{equation}
{}^A{}\widetilde{\nabla}_s s' := [s,s'] + \nabla_{\rho(s')}s \, .\label{Anabla}
\end{equation}
Note that the second term in Eq.~\eqref{Anabla} induced by the ordinary connection $\nabla$ on $A$ renders the left-hand side $C^\infty(M)$-linear in $s$, as it has to be for an $A$-connection. On the other hand, for $A=M\times \g$ with its canonical flat connection $\nabla$, this term vanishes for constant sections and one recognizes the ordinary adjoint action of $\g$ on itself (mimicked by the constant sections of $M\times \g$).

\vskip 3mm\noindent It is certainly to be expected that gauge invariance restricts the fiber metric ${}^A{}g$ in a natural way generalizing ad-invariance and this turns out to be as described above. However, gauge invariance of the new kinetic term for the gauge fields yields also additional constraints on the \emph{previous} data for the scalar fields,\footnote{This is related to the fact that without the kinetic term for the gauge fields, one does not need the Lie algebroid structure and can live with just an anchored bundle. For the kinetic term of the 1-form gauge fields, we need $(A,\nabla)$ to carry a Lie algebroid structure, or at least to have an extension in terms of a Lie${}_\infty$ algebroid.} 
in particular on the Lie algebroid structure of $A$ and the connection $\nabla$ \cite{Mayer-Strobl09,Kotov-Strobl_CYMH15,withyoungpeople}: Denoting by ${}^A{}\widetilde{T} \in \Gamma(A \otimes \Lambda^2 A^*)$ the $A$-torsion of Eq.~\eqref{Anabla}, the required compatibility condition can be cast into the form:
\begin{equation} \label{Szero}
\nabla \left({}^A{}\widetilde{T}\right) =2 \mathrm{Alt} (\iota_\rho   F)  \, ,
\end{equation}
where $F \in \Gamma(A \otimes A^* \otimes \Lambda^2 T^*M)$ denotes the curvature of $\nabla$.

\vskip 3mm\noindent
Despite obtained in a very specific context, the compatibility condition \eqref{Szero}  is natural to consider for \emph{any} Lie algebroid $A$ equipped with a connection---as we will demonstrate also in the body of the paper below from a purely geometrical prespective.  We will thus consider this compatibility of the Lie algebroid structure with the connection also for other geometric structures ${\cG}$ in $(A,\nabla,{\cG})$. 

\vskip 3mm\noindent Eqs.~\eqref{Ag}, \eqref{Anabla}, and \eqref{Szero} were found already quite some time ago \cite{Mayer-Strobl09} in the context of sigma models. However, they came together with another condition, namely that $\nabla$ should be flat, $F=0$. Together with Eq.~\eqref{Szero} this would imply that locally $A$ has to be isomorphic to an action Lie algebroid; indeed, ${}^A{}T$ is a tensorial version of the Lie bracket of $A$ and, if covariantly constant, it implies that the Lie bracket of constant sections is constant. This led the authors of \cite{Mayer-Strobl09} to exclude a further study of this type of theories, since, together with $F=0$ and Eq.~\eqref{Szero}, such theories
would boil down to just the well-known standard Yang-Mills gauge theories (coupled to scalar fields). Only recently \cite{Kotov-Strobl_CYMH15} a way to circumvent this restriction was found, relaxing the condition on $F$: Consider an $A$-valued 2-form $b \in \Gamma(\Lambda^2 T^*M \otimes A)$; its pull-back by $X \colon \Sigma \to M$ can be added to the (covariantized) curvature 2-form of the  gauge field ${\cAfeld}$. Then, gauge invariance of the correspondingly modified kinetic term for ${\cAfeld}$ requires, in addition to Eq.~\eqref{Ag} and Eq.~\eqref{Szero}, the following interesting condition:
\begin{equation} \label{Fb}
F + \left(\iota_\rho \circ \mathrm{D} + \mathrm{D} \circ \iota_\rho \right) b - {}^A{}\widetilde{T}(b,\cdot) = 0 \, .
\end{equation}
The terms involving $b$ are reminiscent of Eq.~\eqref{Lie1}; just since $b$ also has an $A$-component, the covariantized Lie derivative is accompanied by a covariant version of the adjoint action. In particular, for $A=M\times \g$ with its flat connection, this equation implies simply that $b$ is $\g$-invariant, $b \in \Omega^2(M,\g)^\g$. We intend to come back to a geometric meaning of Eq.~\eqref{Fb} elsewhere.



\section{Connections on Lie algebroids}\label{section:connections}

\vskip 3mm\noindent In this section we review the theory of connections on Lie algebroids in the form appropriate for our tasks (in particular, we use the approach to Cartan connections developed in \cite{Blaom04,Blaom05}). We also present this condition as equivalent to the vanishing of a covariant tensor $S$, which is precisely the way it has been found first in the context of gauge theories \cite{Mayer-Strobl09} (but see also \cite{Constrained} for a more recent derivation in the context of the BFV-formalism of constrained systems).

\vskip 2mm\noindent
Every vector bundle $E$ over $M$ gives rise to the short exact sequence \cite{Bott}
\begin{equation} \label{Bott1}
0 \longrightarrow T^*M \otimes E \longrightarrow J^1(E) \longrightarrow E \longrightarrow 0 \, ,
\end{equation}
where $J^1(E)$ is the bundle of first jets of smooth sections of $E$. The embedding of $T^*M \otimes E$ into $J^1(E)$ is determined for every $f,h\in C^\infty (M)$
and $s\in\Gamma (E)$ by the following formula:\footnote{There  is a choice of sign made here, which we find convenient to fix as follows for later purposes.}
 \bqa\label{embegging}
 f\md h\otimes s\mto f \left( hj_1 (s) - j_1 (hs)\right)\,,
 \eqa
where $s\in \Gamma(E)$, $j_1(s)\in\Gamma(J^1(E))$ is the first jet-prolongation of $s$.
Every connection $\nabla$ on $E$ is in one-to-one correspondence with
a splitting $\sigma \colon E \to J^1(E)$ of (\ref{Bott1}): 
 \begin{equation} \label{splitting}
  \sigma(s) = j_1(s) + \nabla s\,,
 \end{equation}
 where  $\nabla s\in \Gamma(T^*M \otimes E)$ is identified with its image in $\Gamma(J^1(E))$.

\vskip 3mm\noindent
Let $(A,\rho,[\cdot,\cdot])$ be a Lie algebroid over $M$. Now (\ref{Bott1}) becomes even an exact sequence of Lie algebroids. The bracket in $J^1(A)$ is defined such that taking the Lie brackets commutes with the prolongation of sections,
 \bqa\label{jets-bracket} [j_1(s), j_1(s')] = j_1([s,s'])
 \eqa
 for all sections $s,s' \in \Gamma(A)$, while its anchor is fixed by the morphism property to obey
 \bqa \label{jet-anchor}
 \rho(j_1(s))=\rho(s)\,.
 \eqa

 \noindent For later purposes, we introduce an $A$-Lie derivative $\cL_s$ along sections $s \in \Gamma(A)$ acting on sections of tensor powers of $A$, $A^*$, $TM$, and $T^*M$ in the following way: $\cL_s s' := [s,s']$,
 $\cL_s X := \cL_{\rho(s)}X \equiv [\rho(s),X]$ for all $s' \in \Gamma(A)$ and $X \in \Gamma(TM)$, extended by the Leibniz rule and the requirement to commute with contractions. Then in particular on
  $T^*M \otimes A$ one has
 \bqa\label{A-Lie}
\cL_{s}(\omega'\otimes s'):=\cL_{\rho(s)}(\omega')\otimes s'+\omega'\otimes [s,s']\,,
\eqa
where $s,s'\in\Gamma (A)$ and $\omega'\in\Omega^1 (M)$.
It follows from (\ref{embegging}), (\ref{jets-bracket}), and (\ref{jet-anchor})  that one has
\bqa\label{cross-bracket}
[j_1 (s), \omega'\otimes s'] &=& \cL_{s}(\omega'\otimes s')\,.
\eqa
 The Lie algebroid structure on $T^*M \otimes A$ is induced by this; it is a bundle of Lie algebras since it belongs to the kernel of the anchor map and
\bqa\label{kernel-bracket}
[\omega \otimes s, \omega' \otimes s'] &=& \langle \omega, \rho(s')\rangle \, \omega' \otimes s - \langle \omega', \rho(s)\rangle \, \omega \otimes s'\,.
\eqa
Note that in the case of the the standard Lie algebroid $TM$, the kernel in \eqref{Bott1}
becomes isomorphic to $\End(TM)$ and the above formula Eq.~\eqref{kernel-bracket} reduces merely to the commutator of endomorphisms.

\vskip 3mm\noindent The anchor map $\rho\colon A\to TM$ admits the first order prolongation $\rho^{(1)}\colon J^1(A)\to J^1 (TM)$, commuting with the jet prolongation of sections,
which makes (\ref{rho-commutes-with-prolongation}) commutative
\bqa\label{rho-commutes-with-prolongation}
\xymatrix{ \Gamma(J^1 (A)) \ar[r]^{\rho^{(1)}}  & \Gamma\left(J^1 (TM)\right)\\
    \Gamma (A)  \ar[u]^{j_1}\ar[r]^{\rho} & \Gamma (TM) \ar[u]^{j_1}   }
\eqa
and which is also a Lie algebroid morphism. Moreover, every Lie algebroid morphism $\varphi\colon (A_1,\rho_1)\to (A_2,\rho_2)$ of Lie algebroids over the same base and the base map being a diffeomorphism admits the first jet prolongation $\varphi^{(1)}\colon J^1(A_1)\to J^1(A_2)$, which is again a Lie algebroid morphism and which commutes with the prolongation of sections, such that the following diagrams are commutative:\footnote{The first diagram is a straightforward generalization of diagram \eqref{rho-commutes-with-prolongation} to which it reduces for $\varphi := \rho \colon A \to TM$. The second diagram results from a likewise diagram for Lie algebroid morphisms over the identity and the fact that the anchor commutes with the prolongation, as is expressed by diagram \eqref{rho-commutes-with-prolongation}.}
\bqa\label{LA-morphism-prolongation}
\xymatrix{ \Gamma(J^1 (A_1)) \ar[r]^{\varphi^{(1)}}  & \Gamma\left(J^1 (A_2)\right)\\
    \Gamma (A_1)  \ar[u]^{j_1}\ar[r]^{\varphi} & \Gamma (A_2) \ar[u]^{j_1}   }
\hspace{15mm}
\xymatrix{ J^1(A_2)  \ar[r]^{\rho_2^{(1)}}  & J^1 (TM)\\
    J^1(A_1) \ar[u]^{\varphi^{(1)}} \ar[r]^{\rho_1^{(1)}}  & J^1 (TM) \ar[u]^{\Id}   }
\eqa


\begin{deff}\label{Cartan-def} $(A,\nabla)$  is called a \emph{Cartan Lie algebroid} over $M$, if $A$ is a Lie algebroid, $\nabla$ a connection on $A\to M$, and its induced splitting $\sigma \colon A \to J^1(A)$ is a Lie algebroid morphism.
 \end{deff}

\noindent The compatibility of a Lie algebroid structure with a connection in the above sense is governed by the vanishing of
the compatibility tensor $S$, the curvature of the splitting (\ref{splitting}), defined by the formula $S (s,s')=[\sigma(s),\sigma (s')]-\sigma\left([s,s']\right)$, where $s,s'\in\Gamma (A)$.
Given that $\rho\left(S (s,s')\right)=0$, it follows from \eqref{Bott1} for $E\sim A$ that $S$ can be identified with a section of $T^*M\otimes A\otimes\Lambda^2A^*$.
The next formula appears in \cite{Blaom04}, Section 2.3 in a slightly different notation.

\begin{proposition} \label{curvature_of_splitting}
\bqa\label{curvature_formula1}
S(s,s')=\cL_{s}\left(\nabla s'\right) - \cL_{s'}\left(\nabla s\right) - \nabla_{\rho(\nabla s)} s' + \nabla_{\rho(\nabla s')}s - \nabla [s,s']   \,.
\eqa
\end{proposition}

\noindent\proof From (\ref{splitting}) we obtain
\beq\nonumber S (s,s')=\left[ j_1 (s)+\nabla s, j_1 (s')+\nabla s' \right] -  j_1 \left( [s, s']\right) -\nabla [s,s']\,.
\eeq
Combining (\ref{jets-bracket}) and (\ref{A-Lie}) we conclude that
\beq\nonumber
S (s,s')=\cL_{s}\left(\nabla s'\right) - \cL_{s'}\left(\nabla s\right) + \left[ \nabla s,\nabla s' \right]- \nabla [s,s'] \,.
\eeq
We finally compute $\left[ \nabla s,\nabla s' \right]$ using (\ref{kernel-bracket}). This accomplishes our proof.
$\square$

\begin{cor}\label{cor-S-in-frame} Let $\{e_a\}$ be a local frame for $A$, such that
$ \rho (e_a)=\rho_a$ and $\nabla e_a=\omega_a^b e_b$,
where $\rho_a$ and $\omega_a^b$ are local vector fields and 1-forms, respectively. Then
\bqa\label{curvature_in_frame}
S (e_a,e_b) &=& \left (\cL_{\rho_a}\left(\omega_b^c\right) - \cL_{\rho_b}\left(\omega_a^c\right) - \left(\iota_{\rho_q}\omega_b^c\right)\omega_a^q
+ \left(\iota_{\rho_q}\omega_a^c\right)\omega_b^q
\right) e_c  \\ \nonumber &+& \omega_b^q [e_a,e_q]-\omega_a^q [e_b,e_q]
 - \nabla [e_a,e_b]   \,.
\eqa
\end{cor}

\noindent\proof By Proposition \ref{curvature_of_splitting}, formula (\ref{curvature_formula1}), we have
\bqa\nonumber
S (e_a,e_b)=\cL_{e_a}\left(\omega_b^c e_c\right) - \cL_{e_b}\left(\omega_a^c e_c\right) - \omega_a^q \nabla_{\rho_q}(e_b)+
\omega_b^q \nabla_{\rho_q}(e_a) - \nabla [e_a,e_b]   \,.
\eqa
From the explicit expression of the A-Lie derivative (\ref{A-Lie}) we immediately obtain (\ref{curvature_in_frame}).
$\square$


\begin{deff}\label{A-connection} An \emph{$A$-connection}  on a vector bundle $V$ is a $C^\infty (M)$-linear map ${}^A{}\nabla$ on sections of $A$
with values in $1$st order differential operators on $V$, which obeys the condition
\bqa\label{A-Leibniz}
{}^A{}\nabla_s (fv)=\rho(s)(f)v+ f{}^A{}\nabla_s (v)
\eqa
for every $s\in\Gamma (A)$ and $v\in\Gamma (V)$.
The \emph{curvature} of an $A$-connection is defined  as
\bqa\label{A-curvature}
{}^A{}R (s,s')= [{}^A{}\nabla_s, {}^A{}\nabla_{s'}]-{}^A{}\nabla_{[s,s']}\,,
\eqa
where $s,s'\in\Gamma (A)$. An $A$-connection is called flat if its $A$-curvature vanishes.
A flat $A$-connection on $V$ gives us a Lie algebroid \emph{representation} of $A$ on $V$.
\end{deff}

\noindent Given a vector bundle $V$, the Atiyah algebroid of $V$, denoted by $L(V)$, is a bundle over $M$,
whose sections are infinitesimal automorphisms of $V$; $L(V)$
is a transitive Lie algebroid, the kernel of the anchor map of which coincides with the bundle of endomorphisms of $V$.
It is easy to verify that an $A$-connection on $V$ is a vector bundle map $A\to L(V)$ which commutes with the
corresponding anchor maps. An $A$-connection is flat if and only if, in addition, this bundle map is a Lie algebroid morphism.

\begin{deff}\label{A-torsion} Given an $A$-connection ${}^A{}\nabla$ on the vector bundle $A$ itself, the $A$-torsion of ${}^A{}\nabla$
is a section of $\Lambda^2 A^*\otimes A\simeq \Hom (\Lambda^2 A, A)$ defined at all $s,s'\in\Gamma (A)$ according to
\bqa\label{A-torsion-formula}
{}^A{}T (s,s')={}^A{}\nabla_s s'- {}^A{}\nabla_{s'}s - [s,s']\,.
\eqa
The dual $A$-connection ${}^A{}\nabla^*$ is determined at all $s,s'\in\Gamma (A)$ by the formula
\bqa\label{A-conn-dual}
{}^A{}\nabla^*_s s'=[s,s']+{}^A{}\nabla_{s'}s\,.
\eqa
\end{deff}
\begin{rem} An easy computation shows that the duality (\ref{A-conn-dual}) is a reflexive operation, i.e. $\left({}^A{}\nabla^*\right)^*={}^A{}\nabla$.
From (\ref{A-torsion-formula}) and (\ref{A-conn-dual}) it follows that the dual $A$-connection has the opposite $A$-torsion; thus,
 ${}^A{}\nabla$ coincides with its dual if and only if ${}^A{}T$ vanishes identically.
\end{rem}

\noindent Let us observe that the Lie derivative of any tensor field $\chi$ along a vector field $X$ depends on the first jet-prolongation of $X$ only, which allows us to introduce
a natural Lie algebroid representation of $J^1 (TM)$ on arbitrary tensor fields, such that
$ j_1 (X)$ acts by the Lie derivative.
The preceding observation, when looked at from a more general
point of view, leads to the representation $\a$ of $J^1(A)$ on the tensor powers of $A$ and $A^*$ for any Lie algebroid $A$, such that for
all $s\in\Gamma(A)$ one has on the one hand
\bqa\label{J1-rep-in-A}
\a\circ j_1(s) = \cL_s\,,
\eqa
and on the other hand, as follows from (\ref{embegging}), for all $s,s'\in\Gamma (A)$, $\omega\in\Omega^1 (M)$
\bqa\label{J1-ker-rep-A}\alpha(\omega\otimes s)s' =\langle \omega, \rho(s')\rangle s\,.\eqa
Correspondingly, $\alpha (\omega\otimes s)$ acts by $- \rho^*(\omega) \vee \iota_s$ on $\mathrm{Sym}^\bullet (A^*)$ and by $- \rho^*(\omega) \wedge \iota_s$  on $\Lambda^\bullet (A^*)$.

\begin{example}\label{J1(TM)} $J^1 (TM)$ is isomorphic to the Atiyah algebroid of $TM$ by means of Eqs.~(\ref{J1-rep-in-A}) and (\ref{J1-ker-rep-A}). The Bott-sequence \eqref{Bott1} specializes here to
\begin{equation} \label{Bott2}
0 \longrightarrow \End(TM) \longrightarrow J^1(TM) \longrightarrow TM \longrightarrow 0 \, .
\end{equation}
 \end{example}

\noindent Combining $\rho^{(1)}\colon J^1(A) \to J^1(TM)$ with the isomorphism from Example \ref{J1(TM)}, we obtain a canonical representation $\tau$ of $J^1(A)$ on the tensor powers of $TM$ and $T^*M$, so that for all $s\in\Gamma (A)$, $X\in\Gamma (TM)$, $\omega\in\Omega^1 (M)$ one has
\bqa\label{J1-rep-in-TM}
  \tau \circ j_1(s) X =[\rho(s),X]\,, \hspace{3mm} \tau(\omega\otimes s) X =\langle \omega, X\rangle \rho(s)\,.
\eqa

\noindent Likewise, $\tau (\omega\otimes s)$ acts by $- \omega \vee \iota_{\rho(s)}$ on $\mathrm{Sym}^\bullet (TM^*)$
and by $- \omega \wedge \iota_{\rho(s)}$ on $\Lambda^\bullet(T^*M)$.

\vskip 3mm\noindent A more general statement is contained in the next proposition:

\begin{proposition}\label{A-conn-from-jet-rep} Let $\mu$ be a representation of $J^1(A)$ on a vector bundle $V$ and let $\nabla$ be a connection on $A\to M$. If we identify $\nabla$ with a splitting $\sigma$ by (\ref{splitting}), then the composition $\mu\circ\sigma$ gives us an $A$-connection on $V$, denoted by ${}^\mu{}\nabla$,
which is flat if and only if the compatibility tensor
$S$ obeys the condition $\mu\circ S(s,s')=0$ for all $s,s'\in\Gamma (A)$. In particular, if $\nabla$ is a Cartan connection, so that $S\equiv 0$, then ${}^\mu{}\nabla$ is flat for every $\mu$. The representations (\ref{J1-rep-in-A}) and (\ref{J1-rep-in-TM}), combined with a connection $\nabla$ on $A$, give us $A$-connections ${}^\a{}\nabla$ and ${}^\tau{}\nabla$ on $A$ and $TM$, respectively, such that for all $s,s'\in\Gamma (A)$, $X\in\Gamma (TM)$ one has
\beq\label{alpha-conn}
 {}^\a{}\nabla_s s' &=& [s,s']+\nabla_{\rho(s')}s \\ \label{tau-conn}
 {}^\tau{}\nabla_s X &=& [\rho(s),X]+\rho\left(\nabla_{X}s\right) \,
\eeq
The anchor map $\rho\colon A\to TM$ obeys the property ${}^\tau{}\nabla\circ \rho=\rho\circ {}^\a{}\nabla$.
\end{proposition}

\noindent\proof It follows from the definition of ${}^\mu{}\nabla$ and the compatibility tensor $S$ along with the explicit formulas (\ref{J1-rep-in-A})(\ref{J1-ker-rep-A})
for the representation $\a$, (\ref{J1-rep-in-TM}) for the representation $\tau$
and (\ref{splitting}) for the splitting $\sigma$. We leave the details to the reader. $\square$


\vskip 3mm\noindent A connection $\nabla$ induces another $A$-connection ${}^A\nabla_\cdot = \nabla_{\rho(\cdot)}$, the $A$-torsion of which will be denoted by ${}^A{}T$.
From (\ref{alpha-conn}) we see that ${}^A\nabla$ is dual to ${}^\a{}\nabla$. Denoting by $F \in \Gamma(\Lambda^2 T^*M \otimes A\otimes A^*)$ the curvature of $\nabla$, we obtain another expression for the compatibility tensor\footnote{It is worth noting that $S$ has appeared in this form in \cite{Mayer-Strobl09}.}.

\begin{proposition}\label{compatibility_vs_S} Let $A$ be a Lie algebroid and $\nabla$ a connection on $A\to M$. Then the compatibility tensor $S$ admits the following expression:
\begin{equation}
S := \nabla \left({}^A{}T\right) + 2\mathrm{Alt} \langle \rho ,  F  \rangle\,.
\end{equation}
Here the anchor $\rho \colon A \to TM$ is viewed as a section of $A^*\otimes TM$, the contraction is taken by the natural pairing $TM \otimes \Lambda^2 T^*M \to T^*M$, $v \otimes \alpha \mapsto \alpha(v, \cdot)$, and the anti-symmetrization is taken over $A^* \otimes A^*$
\end{proposition}
\noindent\proof Our proof starts with the observation that for every $s,s'\in\Gamma (A)$ one has
\beq\nonumber
\nabla \left({}^A{}T \right)(s,s')=\nabla \left( {}^A{}T (s,s')\right)- {}^A{}T \left( \nabla s,s'\right)- {}^A{}T \left( s,\nabla s'\right)\,.
\eeq
The first term of the r.h.s. of the above formula reads as follows
\beqn\nonumber \nabla \left( {}^A{}T (s,s')\right)=
\nabla \left(   \nabla_{\rho (s)}s' -  \nabla_{\rho (s')}s -[s,s'] \right)= \nabla \left(   \nabla_{\rho (s)}s' \right) -
\nabla \left(   \nabla_{\rho (s')}s \right)-\nabla [s,s']\,,
\eeq
while the second and the third terms can be expressed in the form
\beqn\nonumber
- {}^A{}T \left( \nabla s,s'\right) &=& - \nabla_{\rho (\nabla s)}s' +\nabla_{\rho (s')}\left(\nabla s\right) -\cL_{s'}\left(\nabla s\right) \\ \nonumber
- {}^A{}T \left( s,\nabla s'\right) &=& \nabla_{\rho (\nabla s')}s - \nabla_{\rho (s)}\left(\nabla s'\right) + \cL_{s}\left(\nabla s'\right)\,.
\eeq
Here the covariant derivative $\nabla$ is extended to differential forms with values in sections of $A$  by the Leibniz rule. Therefore, combining the above formulas
and using (\ref{curvature_formula1}),
we get
\beq\nonumber
\nabla \left({}^A{}T \right)(s,s')= S(s,s')+ \nabla \left(   \nabla_{\rho (s)}s' \right) - \nabla_{\rho (s)}\left(\nabla s'\right) - \nabla \left(   \nabla_{\rho (s')}s \right)
+ \nabla_{\rho (s')}\left(\nabla s\right)\,.
\eeq
On the other hand, for every vector field $X$ and section $s$ one has
\beqn
\iota_X F (s) =\iota_X \nabla^2 s= \left(\iota_X \nabla + \nabla\iota_X\right)\nabla s-\nabla \left( \iota_X\nabla s\right)=
\nabla_X \left(\nabla s\right) - \nabla\left(\nabla_X s\right)\,.
\eeq
 Finally we obtain
 \beqn
 \nabla \left({}^A{}T \right)(s,s')= S (s,s') -\iota_{\rho (s)} F(s')+ \iota_{\rho (s')} F(s)
 \eeq
or, equivalently,
\beqn
 S(s,s') =\nabla \left({}^A{}T \right)(s,s')+\iota_{\rho (s)} F(s')- \iota_{\rho (s')} F(s) \,.
 \eeq
This accomplishes the proof. $\square$

\begin{cor} \label{CartanLie}
$(A,\nabla)$  is a Cartan Lie algebroid iff $S=0$, i.e.~iff $(A,\nabla)$ satisfies Eq.~\eqref{Szero}.\footnote{We recall that ${}^A{}\widetilde{T} = -{}^A{}T$.}
\end{cor}

\noindent
Let us specify the compatibility of a connection with a Lie algebroid structure for some typical cases.

\begin{example}\label{Action-Cartan-Lie}
Let $A=M\times \g$ be an action Lie algebroid. Then the canonical flat connection $\nabla$ is compatible. Furthermore, \emph{every} flat Cartan Lie algebroid  $(A,\nabla)$ is locally an action Lie algebroid; in fact, one even has:
\end{example}

\begin{proposition}\label{Action-Cartan-Lie_proposition}
Let $(A,\nabla)$ be a flat Cartan Lie algebroid. Then  every point $x \in M$ permits a neighborhood $x \ni U \subset M$ over which there exists an action Lie algebroid $C = U \times \g_U$ such that $(A|_U,\nabla)$ and $(C,\nabla_{canonical})$ are isomorphic as Cartan Lie algebroids. Moreover, within a connected component $M_i$ of $M$, the Lie algebra $\g_U$ does not depend on the choice of $U$,  $\g_U \cong \g_i$.
\end{proposition}

\noindent \proof By Eq.~(\ref{curvature_formula1}) and Corollary \ref{CartanLie}, the Lie bracket of constant sections is constant as well; choosing a covariantly constant frame $e_a$ over some $U \ni x$, $\nabla e_a =0$, the above mentioned property implies that $\md (C^a_{bc}) =0$. This proves the first part, i.e. the local isomorphism of Killing Lie algebroids $A|_U \cong U \times \g_U$. For any other neighborhood $U'$ with constant frame $e_a'$ and non-trivial intersection $U \cap U' \neq \emptyset$, $(e_a)$ and $(e_a')$ are related by an $\R$-linear basis change, which does not change the Lie algebra they generate. $\square$

\begin{example}
If $A$ is a bundle of Lie algebras, i.e.~if $\rho \equiv 0$, then $\nabla$ on $A$ is compatible if and only if it preserves the fiber-wise Lie algebra bracket: $\nabla_X [ \mu, \nu ] = [ \nabla_X \mu , \nu ] + [ \mu, \nabla_X \nu ]$ for all $X \in \Gamma(TM)$ and $\mu,\nu \in \Gamma(A)$.
\end{example}

\begin{example}
If $A=TM$ is the standard Lie algebroid, a connection $\nabla$ on $TM$ is compatible if and only if its dual connection $\nabla^*$ is flat. So, $(TM,\nabla)$ being a Cartan Lie algebroid implies that $M$ is the quotient of a  parallelizable manifold by a properly discontinuously acting discrete group. If, in addition, $\nabla$ is torsion-free, it is self-dual, $\nabla =\nabla^*$, and thus it needs to be flat itself.
\end{example}

\begin{example}[\cite{Blaom05}]
Any torsion-free connection on $TM$ gives rise to a compatible (Cartan) connection on $J^1 (TM)$.
\end{example}

\vskip 3mm\noindent
For what will come in the subsequent section, the following adaptations of Example \ref{J1(TM)} will be of interest: Consider a manifold $M$ equipped with a Riemannian metric $g$. Denote by $L(g)$ the Lie subalgebroid of $J^1 (TM)$ whose sections preserve $g$ by means of the representation of $J^1(TM)$ on tensors.
Into the bargain, $L(g)\subset J^1 (TM)$ is to be viewed as the differential equation whose solutions are Killing vector fields on $(M,g)$, that is, for any section $X\in \Gamma(TM)$, $j_1 (X)\in L(g)$ if and only if $\cL_X (g)=0$. Even though a Riemannian manifold may not possess any Killing vector field, the algebroid $L(g)$ still exists. In this case only there are no first jet prolongations $j_1(X)$
that lands inside the subbundle $L(g)\subset J^1 (TM)$.

\vskip 3mm\noindent
$L(g)$ is naturally isomorphic to
the Lie algebroid of infinitesimal bundle isometries of $(TM,g)$.
Therefore, $L(g)$ is a transitive Lie algebroid, the kernel of the anchor map of which coincides with the bundle of skew-adjoint operators in $TM$,
 canonically identified with $\Lambda^2 (T^*M)$ by use of the metric $g$, i.e.~one has the following short exact sequence of Lie algebroids:
\begin{equation} \label{Killing_equation_algebroid}
0 \longrightarrow \Lambda^2 (T^*M)\longrightarrow L(g) \longrightarrow TM \longrightarrow 0 \, .
\end{equation}

\begin{lemma}\label{Civita_splitting}
The Levi-Civita connection ${}^g{}\nabla$ provides a splitting of (\ref{Killing_equation_algebroid}).
\end{lemma}

\noindent \proof
Let $\sigma_{g}$ be the splitting of the Bott exact sequence (\ref{Bott2}), determined by the Levi-Civita connection as in (\ref{splitting}),
then for all vector fields $X,Y,Z$ on $M$ one has
\beqn
\alpha\circ j_1 (X) (g)  (Y,Z) = \left(\cL_X g \right)(Y,Z)=\cL_X \left( g(Y,Z)\right) - g([X,Y],Z)-g(Y,[X,Z])\,.
\eeq
Given that ${}^g{}\nabla$ is a torsion-free connection, $[X,Y]={}^g{}\nabla_X Y - {}^g{}\nabla_Y X$, and since by formula (\ref{J1-ker-rep-A}) being applied to the standard Lie algebroid $-\a\left({}^g{}\nabla X\right) (g)(Y,Z)=g\left( {}^g{}\nabla_Y X, Z\right)+ g\left( Y, {}^g{}\nabla_Z X\right)$, we obtain
\beqn
\alpha\circ j_1 (X) (g)  (Y,Z) 
={}^g{}\nabla_X (g)(Y,Z)-\a\left({}^g{}\nabla X\right) ( g)(Y,Z)\,.
\eeq
On the other hand, ${}^g{}\nabla (g)=0$,
thus
\beqn\alpha\circ\sigma_{g} (X) (g)=\alpha\left(j_1(X)+{}^g{}\nabla X\right) (g)=0\,.\eeq
Hence it follows that $\sigma_{g}$ takes values in $L(g)$, which is the desired conclusion. $\square$

\section{Lie algebroids over Riemannian manifolds}\label{section:Killing-Lie}


\begin{deff} Let $(A,\rho,[\cdot,\cdot])$  be a Lie algebroid over a Riemannian manifold $(M,g)$ and $\nabla$ a connection on $A$. Then $(A,\nabla)$ and
$(M,g)$ are called \emph{compatible}, if
\begin{equation}\label{Lie-Killing1}
{}^\tau{}\nabla (g)=0 \, ,
\end{equation}
where the $A$-connection ${}^\tau{}\nabla$ is defined as in Eq.~\eqref{tau-conn}. We call the triple $(A,\nabla,g)$ a
\emph{Killing Lie algebroid}, and, if in addition,
also $A$ and $\nabla$ are compatible in the sense of Definition \ref{Cartan-def} (cf.~also Corollary \ref{CartanLie}), a \emph{Killing Cartan Lie algebroid}.
\label{CartanRiemann-def}
 \end{deff}
\noindent \begin{rem}
Eqs.~\eqref{tau-conn} and \eqref{Lie-Killing1} are meaningful also for merely anchored bundles:
\begin{deff}
$(E,\rho)$ is called an \emph{anchored bundle} over $M$, if $E\to M$ is a vector bundle and
$\rho\colon E\to TM$ a vector bundle morphism. We call $(E,\rho,g)$ a \emph{Killing anchored bundle} if Eq.~\eqref{Lie-Killing1} holds true.
\end{deff}
\noindent In some cases of what will follow in this and the subsequent section, we will assume the bundle to carry a Lie algebroid structure, although  sometimes weaker conditions are sufficient for the statement. On the other hand, we will not require the compatibility of this Lie algebroid structure with the connection except if otherwise stated. One option for a passage from a Killing anchored bundle to a Killing Cartan Lie algebroid is explained in \cite{universal}\end{rem}

\begin{lemma}
By means of $g$, we can view the anchor $\rho \in \Gamma(A^* \otimes TM)$ as a section of $A^* \otimes T^*M$, which we denote by $\bar{\rho}$. The connection on $A$ and the Levi-Civita connection on $TM$ induce a connection on  $A^* \otimes T^*M$, which we also denote simply by $\nabla$. Then Eq.~\eqref{Lie-Killing1} holds true if and only if
\begin{equation} \label{rhoKilling}
\mathrm{Sym}\left( \nabla \bar{\rho} \right) = 0 \, .
\end{equation}
\end{lemma}
\noindent\proof One may verify this either by a direct calculation or proceed by the more conceptual consideration that follows:
Using the canonical isomorphism between skew-symmetric bilinear forms
and skew-adjoint operators by means of the Riemannian metric $g$, we deduce that Eq.~(\ref{rhoKilling})
is satisfied if and only if ${}^g{}\nabla \rho (s)-\rho\left(\nabla s\right)$ is a skew-adjoint operator in $TM$ for any $s\in\Gamma(A)$ and, consequently, is a section of $L(g)$.
From (\ref{splitting}) and (\ref{rho-commutes-with-prolongation})
we have
\beq\label{splitting_intermediary}
 \sigma_{g}\circ\rho (s)-\rho^{(1)}\circ\sigma (s)
= {}^g{}\nabla \rho (s)-\rho\left(\nabla s\right)\eeq
since $\rho^{(1)}(\nabla s)=\rho(\nabla s)$, where $\rho$ is extended to act on sections of  $T^*M\otimes A$ as $\mathrm{id}\otimes\rho$.
Thus the l.h.s.~of Eq.~(\ref{splitting_intermediary}) is also a section of $L(g)$.
By Lemma \ref{Civita_splitting},
the image of  $\sigma_{g}$ is contained in $L(g)$, therefore Eq.~(\ref{rhoKilling}) is seen to be equivalent to the requirement  $\rho^{(1)}\circ\sigma \colon A\to L(g)$,
which in turn is equivalent to ${}^\tau{}\nabla$ annihilating the metric $g$; the latter fact follows from the construction of ${}^\tau{}\nabla$ and $L(g)$. $\square$

\begin{lemma} \label{Killing-lemma}
Suppose there exists a section $s\in \Gamma(A)$ which is covariantly constant, $\nabla s = 0$. Then $v := \rho(s)$, provided non-zero, is a Killing vector field of the metric $g$, $\cL_v g = 0$.
\end{lemma}
\noindent \proof Clearly, Eq.~\eqref{rhoKilling} implies
 $\mathrm{Sym} \left(\nabla \bar{v}\right) = 0$ where $\bar{v}\equiv \langle s,\bar{\rho}\rangle \in \Gamma(T^*M)$. Rewriting this in terms of $v$, $\bar{v}=g(v,\cdot)$, this equation is known to become $\cL_v g = 0$. $ \square$


\begin{rem} Using notations from Corollary \ref{cor-S-in-frame}, we can rewrite (\ref{Lie-Killing1}) as \begin{equation} \label{extKill}
\cL_{\rho_a} g-\omega_a^b\vee\iota_{\rho_b}g=0
\end{equation}
for all $a=1, \ldots, \mathrm{rk} A$; this expression coincides with the original form
of the extended Killing equations found in \cite{Kotov-Strobl14}. This fact, as well as the observation in the previous Lemma, led us to use the term ``Killing Lie algebroid''. While this nomenclature reflects well the relation to symmetries of the geometric structure on the base here, its generalization to other appropriately compatible geometric structures on the base (cf.~the subsequent section) is not obvious (for example if the Killing vectors are replaced by  symplectic ones).
\end{rem}
\noindent Since equations such as \eqref{extKill} are in the
form in which they appear in the context of their original appearance, and henceforth we will turn to a reformulation such as in
Eq.~\eqref{Lie-Killing1}, using the $A$-connection ${}^\tau{}\nabla$ also in more general context, we consider it illustrative to show equivalence of these two formulations by an explicit calculation:

\vspace{1mm}\noindent {\bf Proof (}equivalence of \eqref{Lie-Killing1} with the local expression \eqref{extKill}{\bf ).} Let $X\in \Gamma(TM)$ be a vector field on $M$ and $s$ a section in $A$.
${}^\tau{}\nabla_s$ acting on the function  $g(X,X)$ agrees, by definition of an $A$-Lie derivative, with $\cL_{\rho(s)}$ acting on it. On the other hand, evidently we have  ${}^\tau{}\nabla_s\left(g(X,X)\right)= ({}^\tau{}\nabla_s g)(X,X) + 2 g({}^\tau{}\nabla_s X,X)$. Since by the defining equation \eqref{tau-conn}, ${}^\tau{}\nabla_s X = \cL_{\rho(s)} X + \rho(\nabla_X s)$, we obtain $({}^\tau{}\nabla_s g)(X,X) = (\cL_{\rho(s)} g)(X,X) -2g(\rho(\nabla_Xs),X)$. Using $\rho(\nabla_X e_a)= (\iota_X\omega^b_a) \, \rho(e_a)$ for any local frame $e_a$ in $A$ then yields Eq.~\eqref{extKill}.
 $\square$


\begin{example}\label{example_flat_CartanKilling}
Let $\g$ be the isometry Lie algebra of a metric $g$ on $M$. Then the action Lie algebroid $M \times \g$, its canonical flat connection, together with $g$ form a Killing Cartan Lie algebroid. Furthermore, as remarked in Example \ref{Action-Cartan-Lie} and Proposition \ref{Action-Cartan-Lie_proposition}, every flat Cartan Lie algebroid  $(A,\nabla)$ is locally an action Lie algebroid, $(A|_U,\nabla) \cong (U \times \g_i,\nabla_{canonical})$, where the Lie algebra $\g_i$ is the same for any connected component $M_i$ of $M$. Adding the prefix ``Killing'' then implies that now the image $\h_i$ of $\g_i$ by the anchor map $\rho$, a morphism of Lie brackets, is necessarily a sub-Lie algebra  of the local isometry Lie algebra $\mathrm{iso}_i$ of $(M_i,g)$: $\rho(\g_i) = \h_i \subset \mathrm{iso}_i$. The Lie algebras $\g_i$ of the local action Lie algebroids are thus extensions of local isometry Lie algebras: $ 0 \to \ker(\rho) \to \g_i \stackrel{\rho}{\longrightarrow} \h_i \to 0$ in this case.
\end{example}


\begin{example}
The standard Lie algebroid $TM$ of a Riemannian manifold $(M,g)$ is a Killing Lie algebroid with respect to any metrical connection $\nabla$ (e.g.~w.r.t.~the Levi-Civita connection ${}^g{}\nabla$). Indeed, in this case $\bar{\rho} = g$ and Equation \eqref{rhoKilling} holds true even without the symmetrization. $(TM,{}^g{}\nabla,g)$ is a Killing Cartan Lie algebroid iff $g$ is a flat metric, i.e.~iff $R_{{}^g{}\nabla}=0$. $(TM,\nabla,g)$ is a Killing Cartan Lie algebroid iff $\nabla$ is metrical and its curvature $R$ and torsion $T$ satisfy the equation  $T^i_{jk;l} = R^i_{jkl} - R^i_{kjl}$, written for clarity in index notation, the semi-colon denoting the covariant derivative.
\end{example}

\begin{example}[\cite{Kotov-Strobl14}]\label{reg_foliation} A characteristic example results from (regular) foliations $\cF$ of a manifold $M$ with a smooth quotient $Q=M/\cF$. There \emph{exists} a connection $\nabla$ on $A=T\cF$ such that
$(T\cF,\nabla,g)$ forms a Killing Lie algebroid, \emph{iff} $(M,\cF,g)$ provides a Riemannian submersion, i.e.~\emph{iff} $Q$ can be equipped with a metric $g_Q$ such that for any point $p \in M$ the natural projection $(T_p\cF)^\perp \to T_{[p]}Q$ is an isometry.
\end{example}

\noindent While on the one hand Killing Lie algebroids (and their generalizations to other geometric structures) encompass symmetries in terms of a connection such that locally flat sections (with non-trivial image w.r.t.~the anchor) reproduce infinitesimal, local symmetries, on the other hand, even in the absence of symmetries, \emph{if} the leaf space $Q:=M/\cF$ is smooth, the conditions ensure that there is some kind of quotient construction (of the same kind in the Riemannian case, i.e.~$(Q,g_Q)$ is the quotient of $(M,g)$ by the foliation $\cF$). This aspect is in its spirit closely related to the gauge theories where the notion arose from, cf.~\cite{Kotov-Strobl14} as well as the Introduction; and in both cases, the geometrical as well as the gauge theoretical context, the notion
becomes particularly interesting when there is no smooth quotient, providing a smooth description of it. Evidently, the connection contains information about the metric on the leaves of the foliation factored out in the previous example, cf., e.g., Eq.~\eqref{extKill}. It permits a partial reconstruction of the metric on the total space. More precisely, formulated for a regular foliation with smooth quotient, this looks as follows:
\begin{proposition}
Let $\pi \colon M \to Q$ be a bundle over a Riemannian manifold $(Q,g_Q)$ with connected fibers, let $\nabla_M$ be an Ehresmann connection on $\pi$, i.e.~a splitting $T_xM = H_x \oplus V_x$ with $V_x \equiv T_x\cF$, and $\sigma \in \Gamma(\pi)$ a section such that for all $x$ in its image, $x \in \sigma(Q)$, one is given a smoothly varying metric $\eta_x$ on the vertical subspace $V_x$. Denote by $T\cF\subset TM$ the foliation Lie algebroid over $M$ corresponding to $\pi$. For every choice of a $T\cF-$connection
(a connection along the fibers) ${}^{\sst T\cF}\nabla$ on $T\cF$,
 the holonomy of which at $T\cF|_{\sigma (Q)}$ preserves the metric $\eta$,
 there is a unique connection $\nabla$ on $T\cF$ and a unique metric $g$ on $M$, such that $(T\cF,\nabla,g)$ is a Killing Lie algebroid, $H_x$ is orthogonal to $V_x$ w.r.t.~$g_x$ for all $x \in M$, the restriction of $g_x$ on $V_{\sigma(q)}$ agrees with $\eta_{\sigma(q)}$ for each  $q \in Q$, and the restriction of $\nabla$ on fibers coincides with ${}^{\sst T\cF}\nabla^*$,
 the dual $T\cF-$connection defined by Eq.~\eqref{A-conn-dual}.
\end{proposition}
\noindent \proof The connection $\nabla$ is defined by the requirement that $\nabla_X (s')=[X,s']^V$, $\nabla_s s'={}^{\sst T\cF}\nabla^*_s s'$ for every $X\in\Gamma (H)$ and $s,s'\in\Gamma (T\cF)$, where $[-,-]^V$ is the vertical component of the corresponding Lie bracket.
Use Eq.~\eqref{Lie-Killing1} and the holonomy property, noting that ${}^\tau \nabla$ is a partially defined covariant derivative, permitting one to transport $\eta_{\sigma(q)}$ to the metric $g|_{T\cF}$.
By Eq.~\eqref{tau-conn}, for every connection $\nabla'$ on $T\cF\to M$, the vertical subbundle $T\cF$ is ${}^\tau \nabla'-$invariant.
Moreover, the induced connection ${}^\tau\nabla'/T\cF$ on the normal bundle $TM/T\cF$ gives rise to the canonical normal transversal action of $T\cF$,
thus it is flat and it does not depend on the choice of $\nabla'$.
On the other hand, by the construction of $\nabla$, both subbundles $H$ and $V=T\cF$ are ${}^\tau \nabla-$invariant, so that $(H,{}^\tau \nabla|_H$) is canonically isomorphic to the
transversal representation of $T\cF$. This allows us to define a unique ${}^\tau \nabla-$invariant metric on $H$, such that
$\pi_x\colon H_x \to T_{\pi(x)}Q$ is an isometry for all $x\in M$. The metrics on $H$ and $T\cF$ determine a unique metric on $M$, such that $H$ is orthogonal to $V$ (and thus
$\pi\colon M\to Q$ is a Riemannian submersion). Since the obtained metric on $M$ is ${}^\tau \nabla-$invariant, $(T\cF,\nabla,g)$ is a Killing Lie algebroid.
$\square$

\begin{example} That one may encounter simple obstructions in the search for further examples of Killing Lie algebroids or even Killing anchored bundles is illustrated by the following:
\begin{proposition} \label{rho0}
Let $E=\R^2 \times \R$ with $\rho(x,y,u) = ux^n \frac{\partial}{\partial_x}$ for some $n \in \N$.
There exists \emph{no} connection $\nabla$ on  $E$ and metric $g$ on $M=\R^2$ such that $(E,g,\rho)$ forms a Killing anchored algebroid.
\end{proposition}
\begin{rem}  Note that every anchored line bundle $L \to M$ is compatible with one and only one Lie algebroid structure, i.e.~the anchor $\rho \colon L \to TM$ determines the Lie algebroid structure already completely and there exists one for every choice of the map $\rho$. Indeed, consider a local chart $U \subset M$ over which $L|_U \cong U \times \R$. The unit section $1$ of this trivial bundle is mapped to a vector field $v=\rho(1) \in \Gamma(TU)$. Any other local section $s \in \Gamma(L|_U)$ then satisfies $s = f \cdot 1$ for some $f\in C^\infty(U)$. For two such sections $s,s'$, the bracket takes necessarily the form:
$[s,s']=\left(fv(f')-f'v(f)\right) \cdot 1$.
Since, on the other hand, one verifies that this satisfies all axioms of a Lie algebroid over $U$, and that the bracket is equivariant with respect to a change of basis over $U$, one may extend it also to all of $M$ by a partition of unity; this also proves existence.
\end{rem}
\noindent {\bf{Proof} \rm{(of Proposition \ref{rho0})}\bf{.}} For the constant section $s=1$ with $v:=\rho(1)=x^n \partial_x$ and $\nabla(1) = \omega \otimes 1$, Equation \eqref{Lie1} becomes
\begin{equation}
\cL_v g = \omega \vee \iota_v g\,\; ,
\end{equation}
where $\omega \vee \iota_v g \equiv \omega \otimes \iota_v g+\iota_v g\otimes\omega$. Comparing the components of $\md x \otimes \md x$ on both sides for the given vector field $v$, yields the equation
$$ x^n \partial_x(g_{xx}) + 2 g_{xx}n x^{n-1} = 2 \omega_x g_{xx} x^n $$ or, since $g_{xx} > 0$, equivalently
$$ n = x \left(\omega_x - \frac{1}{2} \partial_x \ln(g_{xx}) \right) \, .$$ For whatever the choice of $g$ and $\omega$, i.e.~$\nabla$, this last equation yields a contradiction upon evaluation at $x=0$.
 $\square$
\end{example}

\noindent All the more it is important to provide further, non-trivial examples of Killing Lie algebroids. That a very large class of them exists follows from the following theorem.
\begin{theorem} \label{proper} Any Lie algebroid $A$ which is integrable to a proper Lie groupoid permits a metric $g$ and a connection $\nabla$ to turn $(A,\nabla,g)$ into a Killing Lie algebroid.
\end{theorem}

\noindent This Theorem is corollary of the result of M.~de Hoyo and R.~Fernandes \cite{delHoyo-Fernandes15}, which states that any proper Lie groupoid admits what they call a 2-metric,
as well as Lemma \ref{lemma_submersion1} and Lemma \ref{lemma_submersion2} below;
the second of those lemmas asserts that the Lie algebroid of any Lie groupoid with 1-metric admits a canonical Killing Lie algebroid structure.\footnote{The latter fact was found first by  Camille Laurent-Gengoux and Sylvain Lavau (communicating it to us without showing us their proof).}

\vskip 3mm\noindent Let $\pi\colon (\tilde M,\tilde g)\to (M,g)$ be a Riemannian submersion. Denote by $V(\pi)$ the subbundle
of $\pi-$vertical vectors and by $H$ its orthogonal complement with respect to ${\tilde g}$, such that every vector field $\tilde X$ on $\tilde M$ admits
the canonical decomposition into the horizontal and vertical parts, $\tilde{X}_H$ and $\tilde{X}_V$, respectively. Denote by $X^h$
the unique horizontal lift of a base vector field $X$, such that $\md\pi (X^h)=X$.
The following 3-tensor field $O$, introduced by O'Neill
in \cite{O'Neill66},\footnote{In \cite{O'Neill66}, this tensor was called $A$;  we refrained from this notation here, since $A$ already denotes the Lie algebroid throughout this paper.} associates to a pair of vector fields $\tilde X$, $\tilde Y$ on $\tilde M$ a vector field $O_{\tilde X} \tilde Y$, where, by definition,
\beq\label{O'Neill_A-tensor}
O_{\tilde X} \tilde Y =\left( {}^{\tilde g}{}\nabla_{\tilde{X}_H} \tilde{Y}_H \right)_V + \left( {}^{\tilde g}{}\nabla_{\tilde{X}_H} \tilde{Y}_V \right)_H\,.
\eeq

\begin{lemma}[O'Neill, \cite{O'Neill66}]\label{lemma_submersion1}{\mbox{}\vskip 0.5mm}
\noindent
1. At each point, $O_{\tilde X}$ is a skew-symmetric linear operator on the tangent space of $\tilde M$ interchanging the horizontal
and vertical subspace.
\vskip 1mm
\noindent 2. If $\tilde X$ and $\tilde Y$ are horizontal vector fields on $\tilde M$, then
$O_{\tilde X} \tilde Y =\frac{1}{2} [\tilde X,\tilde Y]_V$.
\vskip 1mm
\noindent 3. For $H-$horizontal lifts $X^h$, $Y^h$ of basic vector fields $X$, $Y$:
$\left({}^{\tilde g}{}\nabla_{X^h}Y^h\right)_H=\left({}^g{}\nabla_{X}Y \right)^h$.
\end{lemma}

\begin{cor}\label{corollary_submersion1}
If $\tilde X$, $\tilde Y$ are horizontal vector fields on $\tilde M$
and $\xi$ a vertical vector field, then $\tilde g ({}^{\tilde g}{}\nabla_{\tilde X} \tilde Y, \xi )$, which is evidently equal to
$-{\tilde g} ({}^{\tilde g}{}\nabla_{\tilde X}\xi,\tilde Y)$,
is skew-symmetric with respect to $\tilde X$ and $\tilde Y$.
If $X^h$, $Y^h$, and $Z^h$ are $H-$horizontal lifts of vector fields $X$, $Y$, and $Z$ on $M$, respectively, then
\bqa\label{submersion_formula2}
{\tilde g} ({}^{\tilde g}{}\nabla_{X^h}Y^h, Z^h )=g ({}^g{}\nabla_{X}Y, Z )\,.
\eqa
Finally, for $X$, $Y$ vector fields on $M$ and $X^h$, $\tilde Y$ vector fields on $\tilde M$, such that $X^h$ is the horizontal lift of $X$ and $\tilde Y$ is any $\pi-$projectible lift of $Y$,
i.e.~a vector field on $\tilde M$ such that $\md\pi (\tilde Y)=Y$, one has
\bqa\label{submersion_formula3}
{\tilde g} \left({}^{\tilde g}{}\nabla_{X^h} \tilde Y, X^h\right)=g \left({}^g{}\nabla_{X} Y, X\right)\,.
\eqa
\end{cor}
\noindent \proof To see the skew-symmetry, use the defining equation for the  O'Neill's tensor $O$, Eq.~\eqref{O'Neill_A-tensor}, and the second item in O'Neill's Lemma \ref{lemma_submersion1}. Similarly, Eq.~\eqref{submersion_formula2} follows from the third part of the Lemma and the definition of a Riemannian submersion. The final equation, Eq.~\eqref{submersion_formula3}, then results from Eq.~\eqref{submersion_formula2} and the before-mentioned skew-symmetry.$\square$

\begin{lemma}\label{lemma_submersion2} Let $\cG$ be a Lie groupoid over $M$ with the source map $s$, target map $t$, and the identity bisection $e$.
Assume that $\cG$ is endowed with a Riemannian metric $\eta$, such that both $s$ and $t$ are Riemannian submersions $s\colon (\cG,\eta)\to (M,g)$ and
$t\colon (\cG,\eta)\to (M,g')$ for some metrics $g$ and $g'$ on $M$, respectively.\footnote{The base metrics $g$ and $g'$ do not need to coincide. This condition is weaker than a $1-$metric on $\cG$
in the original notations from \cite{delHoyo-Fernandes15}, where $g$ and $g'$ must coincide.
Thus we essentially prove a stronger statement. A proper Lie groupoid admits even a $2$-metric which contains far more information than what we need for the existence of a Killing Lie structure; we shall develop this subject in the next paper.}
Then the Lie algebroid $A$ of $\cG$, identified with
left invariant $t$-vertical vector fields on $\cG$, is a Killing Lie algebroid over $(M,g)$.
\end{lemma}
\noindent\proof Let $\xi$ be a section of $A$, identified with the corresponding vector field on $\cG$, such that $\rho(\xi)=\md s(\xi)$, $X$ be a vector field on $M$,
and $X^h$ be the horizontal $s$-lift of $X$, orthogonal to the $s$-fibers. From (\ref{submersion_formula3}) we obtain
\bqa\label{groupoid_equality1}
\eta \left({}^\eta{}\nabla_{X^h} \xi, X^h\right)=g \left({}^g{}\nabla_{X} \rho(\xi), X\right)\,.
\eqa
Consider now the following symmetric $2$-form $C(\tilde X,\tilde Y)= \eta({}^\eta{}\nabla_{\tilde X} \xi, \tilde Y)$ on $\cG$.
Using metric $\eta$, we identify $C$ with a section of $T\cG\otimes T\cG$;
the latter decomposes into the direct sum of orthogonal components according to the orthogonal decomposition $T\cG=V(t)\oplus V(t)^\perp$, where $V(t)=\Ker \,\md t$.
By the the first part of Corollary \ref{corollary_submersion1} being applied to the Riemannian submersion $t$,
we get that $C=C_0+C_1+ C_2$, where $C_0$, $C_1$ and $C_2$  are sections of $\Lambda^2 V(t)^\perp$,
$T\cG\otimes V(t)$, and $V(t)\otimes V(t)^\perp$, respectively.
Obviously, $C_1 (\tilde X,\tilde Y)=\eta(\nabla^{pr}_{\tilde X} \xi, \tilde Y)$,
where $\nabla^{pr}$ is the orthogonal projection of the Levi-Civita connection onto sections of $V(t)$.
Besides, there exists a unique section $\psi^{pr}$ of the bundle of
endomorphisms $T\cG\to \End \left( V(t)\right)$, which factors through the orthogonal projection onto $V(t)^\perp$, such that
$C_2(\tilde X,\tilde Y)=\eta (\tilde X,\psi^{pr}_{\tilde Y} \xi)$.
Since $C_0$ is a skew-symmetric 2-form,
we conclude that $C(\tilde X,\tilde X)=\eta (\nabla^{tot}_{\tilde X}\xi,\tilde X)$,
where $\nabla^{tot} =\nabla^{pr}+\psi^{pr}$. Now (\ref{groupoid_equality1}) reads as follows:
\bqa\label{groupoid_equality2}
g \left({}^g{}\nabla_{X} \rho(\xi), X\right)=\eta \left(\nabla^{tot}_{X^h}\xi, X^h\right)\,.
\eqa
Given that the l.h.s. of (\ref{groupoid_equality2}) is constant along $s$-fibers, we can evaluate the r.h.s. at the identity bisection $e$. Thus
\beqn
g \left({}^g{}\nabla_{X} \rho(\xi), X\right)=\eta \left(\left(\nabla^{tot}_{ X^h}\xi\right)|_e, X^h|_e\right)=g\left (\rho\circ \nabla_X\xi, X \right)\,,
\eeq
where $\nabla_X$ is defined on sections of $A$ as the composition of $\nabla^{tot}_{ X^h}$ and the evaluation at the identity bisection $e$.
 In this way, we have obtained a connection on $A$ such that
the identity (\ref{rhoKilling}) holds. This concludes the proof of Lemma \ref{lemma_submersion2} and, at the same time, of Theorem \ref{proper}.
$\square$

\begin{rem}
Note that Theorem \ref{proper} does not give any information about an eventual compatibility of the Lie algebroid structure on $A$ with the connection $\nabla$. On the other hand, using averaging methods by fiber-integration developed on proper groupoids, the Theorem does not only have an existence part, it also permits one to construct non-trivial examples of Killing Lie algebroids. We intend to come back to such examples elsewhere.
\end{rem}


\vskip 3mm \noindent We now turn to some properties one may show to hold for Killing Lie algebroids. For example, we claimed before that Example \ref{reg_foliation} is in some sense characteristic, although the regular foliation and the smooth quotient seem very restrictive for a Lie algebroid. A generalization of Riemannian submersions to the non-smooth setting and for also singular foliations is given by Riemannian foliations \cite{Molino88}. It is thus comforting to find
\begin{proposition}\label{Riemannian_foliation} The (possibly singular) foliation of the base $M$ defined by a Killing Lie algebroid $(A,\nabla,g)$ over $M$
becomes a Riemannian foliation relative to $g$.
\end{proposition}
\noindent \proof  Let $p$ be a point of $M$ and let $\gamma$ be a geodesic curve with the natural parameter $s\in [0,1]$ such that $\dot{\gamma} (0)\in (T_p\cF)^\perp$,
where $\dot{\gamma}$ is the derivation along $s$.
Consider the pullback bundles $\gamma^*(A)$ and $\gamma^* TM$ together with the corresponding pullback connections and the pullback of the anchor map regarded as a bundle morphism $\rho\colon \gamma^* A\to \gamma^* TM$.
Given that $\gamma$ is a $1-$dimensional manifold, the pullback connection is flat, thus we can choose a flat trivialization $\{e_a\}$ of $\gamma^*(A)$
in some neighborhood of $p\in\gamma$.
Now we have:
\beq\label{natural_derivation}
 \pt_s g(\dot{\gamma}, \rho(e_a))=g({}^g{}\nabla_{\dot{\gamma}}\dot{\gamma}, \rho(e_a))+g(\dot{\gamma},{}^g{}\nabla_{\dot{\gamma}}\rho(e_a) )\,.
\eeq
The first term of (\ref{natural_derivation}) is identically zero since $\gamma$ is a geodesic curve, while the second term vanishes because of the extended Killing equation, Eq.~(\ref{rhoKilling}) being applied to
the flat frame $\{e_a\}$. Indeed, ${}^g{}\nabla_{\dot{\gamma}}\bar\rho(e_a) = \left(\nabla_{\dot{\gamma}} \bar{\rho}\right) (e_a)$
and thus $g(\dot{\gamma},{}^g{}\nabla_{\dot{\gamma}}\rho(e_a) )=\mathrm{Sym}\left( \left(\nabla \bar{\rho}\right) (e_a)\right) (\dot{\gamma}, \dot{\gamma})=0$.
 Therefore $g(\dot{\gamma}, \rho(e_a))$ does not depend on $s$, thus it must be zero as $\dot{\gamma} (0)\in (T_p\cF)^\perp$. This proves that
the geodesic remains orthogonal to the foliation for all $s$, which is a possible characterization of a Riemannian foliation, cf.~\cite{Molino88}.
$\square$

\vskip3mm \noindent Let us for the rest of the section consider a fixed Riemannian base $(M,g)$ and determine some conditions under which this can be the base of a Killing Lie algebroid. Remember that the Lie algebroid $L(g)$, fitting into the sequence \eqref{Killing_equation_algebroid}, exists for any metric $g$. We will use below the germ version of the pre-image of a vector sub-bundle by a bundle map: namely, given a vector bundle map $\phi\colon V\to W$
between two vector bundles over $M$ and a vector subbundle $W'\subset W$, we define $\phi^{-1}(W')$ at $x\in M$ as the set of vectors $v\in V_x$ which admit prolongations to local sections of $V$ with the image in $W'$.

\begin{proposition}
Let $(A,\rho,[\cdot, \cdot])$ be a Lie algebroid  over a Riemannian manifold $(M,g)$. Then $A$ admits a connection which satisfies the equation (\ref{Lie-Killing1}),
if and only if the preimage of $L(g)$ by $\rho^{(1)}$ is surjective over $A$. The choice of such a compatible connection $\nabla$ is in one-to-one correspondence
with the choice of a splitting of $\left(\rho^{(1)}\right)^{-1}(L(g))\to A$.
\end{proposition}

\noindent\proof  By the above (germ) definition,
the pre-image of $L(g)$ by $\rho^{(1)}$ is surjective over $A$ if and only if for any point $x$ in $M$ there exists an open neighborhood $U$ of $x$
and a local connection on $A|_{U}$, the restriction of $A$ to $U$, which is compatible with $g$. Let us choose an open cover of $M$ together with a local connection on $A$ over each open subset from the cover, which is compatible with $g$. Using a partition of unity subordinated to the chosen open cover and taking the corresponding linear combination of the local connections, we construct a global connection on $A$ compatible with $g$ on $M$. 
$\square$

\vskip 3mm
\noindent Given existence of the Killing Lie algebroid $(A,\nabla)$  over $(M,g)$, one may ask for the ambiguity in the choice for $\nabla$.
Let us for this purpose consider the following Koszul complex of vector bundles:
$\widetilde{C}^q := \mathrm{Sym} T^*M\otimes \Lambda^q A$,  $q=0,\ldots ,\mathrm{rk} A$, with the differential $\tilde\de$
obtained by the natural extension
of $\bar\rho$, where the latter is regarded as a section of $A^*\otimes T^*M\simeq \Hom (A, T^*M)$. More precisely, $\tilde\de$ acts on $\Lambda^q A$ by contraction with the first factor
of $\bar\rho$ and on $\mathrm{Sym}^p T^*M$ by symmetric multiplication on the second one. Now we take the twisted
complex $C^\bullet =\widetilde{C}^\bullet\otimes A^*$ with the differential $\de =\tilde\de\otimes\mathrm{id}$.
This complex is graded by sub-complexes
 $C_k^\bullet =\bigoplus_{p+q=k} \mathrm{Sym}^p T^*M\otimes \Lambda^q A\otimes A^*$, $k=0, \ldots, \dim M+\mathrm{rk} A$. Let $\cC^q_k:=\Gamma \left( C^q_k\right)$
 with the induced differential $\delta\colon \cC^q_k\to \cC^{q-1}_k$ (we identify the bundle map $\de$ with the underlying operator acting on sections).

\begin{proposition} Let $(A,\nabla,g)$ be a Killing Lie algebroid. Then $\nabla+\psi$ for some  $\psi\in\Gamma (T^*M\otimes \End A)$ defines another connection compatible with $g$, if and only if $\delta (\psi)=0$.
\end{proposition}
\noindent \proof
Since $\mathrm{Sym} ( \nabla\bar\rho)=0$, the
 equation (\ref{rhoKilling}) for $\nabla +\psi$ implies that
 \beqn
 \mathrm{Sym}\otimes\mathrm{id} \left((\mathrm{id}\otimes\bar\rho\otimes\mathrm{id})(\psi)\right) =0\,,
 \eeq
 where $\bar\rho$ is regarded
as a section of $\Hom (A, T^*M)$ and $\psi$ as a section of $T^*M\otimes A\otimes A^*$, so that $\bar\rho$ is applied to the second factor of $\psi$.
Let us think of $\psi$ as an element of $\cC^1_2$. By the definition of the Koszul differential $\de$
the following diagram is commutative
\bqa\nonumber
\xymatrix{
TM^*\otimes A\otimes A^* \ar[rr]^{\mathrm{id}\otimes\bar\rho\otimes\mathrm{id}} \ar[drr]_{\de} && T^*M\otimes T^*M \otimes A^* \ar[d]^{\mathrm{Sym}\otimes\mathrm{id}} \\
    && \mathrm{Sym}^2T^*M \otimes A^*
    }
\eqa
This proves the proposition.
$\square $

\vskip3mm \noindent We finally return to the aspect of the relation of Killing Lie algebroids to Killing vectors. It is well-known that the space of Killing vectors on an $n$-dimensional, connected manifold $M$ is a vector space of dimension at most $n(n+1)/2$. There is a straightforward generalization of this fact: Starting from Eq.~\eqref{rhoKilling}, one derives an equation expressing the two-fold covariant derivatives of $\bar{\rho}\equiv \iota_\rho g$ in terms of itself and its first covariant derivative in the standard way; the only difference now is that the curvature appearing in this equation now is the total curvature $F_\nabla + R_{{}^g\nabla}$. Since in this context we never refer to the Lie algebroid structure on the bundle, one obtains a statement for Killing anchored bundles: 
\begin{proposition}
Let $E$ be a rank $r$ vector bundle over an $n$-dimensional, connected Riemannian manifold $(M,g)$ and $\nabla$ a connection on $E$. Denote by $V_g$ the vector space of Killing vector fields on $(M,g)$. Then the following facts hold true: Anchors $\rho \in \Gamma(A^* \otimes TM)$ satisfying \eqref{barrho} form a finite-dimensional vector space $W$, whose dimension is bounded by  $\mathrm{dim} W \leq \frac{r \cdot n(n+1)}{2}$.
For a metric $g$ of constant curvature and $A = M \times \R^r$ with its natural flat connection, $W \cong \left(V_g\right)^{\otimes r}$, attaining the above bound. In general, the dimension of $W$ can be greater than $r$ times the dimension of $V_g$.
\end{proposition}

\noindent The last sentence is proven, e.g., by Example \ref{reg_foliation} with $M=Q \times \cF$ for a metric $g = g_Q+g_{\cF}$ with no isometries: $\dim V_g = 0$ and $\dim W \ge 1$.

\vskip 3mm
\noindent While flat Killing Cartan Lie algebroids are locally action Lie algebroids with $\rho$ being a (Cartan) Lie algebroid morphism into $M \times \mathrm{iso(g)}$, where $\mathrm{iso(g)}$ denotes the isometry Lie algebra of the metric $g$, cf.~Example \ref{example_flat_CartanKilling}, the relation of
the more general flat Killing Lie algebroids to isometries of $g$ is more intricate. We defer some further aspects of flat Killing Lie algebroids, Cartan and not Cartan, to Appendix A.

\section{Lie algebroids over manifolds with other geometric structures}\label{section:symplectic}

\begin{deff}\label{Bilinear-Lie-Killing_def} Let $(A,\rho,[\cdot,\cdot])$  be a Lie algebroid over $(M, \Phi)$, where $\Phi \in \Gamma(T^*M \otimes T^*M)$ is a bilinear form on $M$,
$\nabla$ be a connection on $A$, and $\psi$ a section of $T^*M\otimes \End (A)$.
Then $(A,\nabla,\psi)$ and $(M,\Phi)$ are called \emph{compatible}, if
\begin{equation}\label{Bilinear-Lie-Killing}
{}^\tau{}\nabla^{comb} (\Phi)=0\,,
\end{equation}
where ${}^\tau{}\nabla^{comb}={}^\tau{}\nabla^+\otimes\Id +\Id\otimes {}^\tau{}\nabla^-$ corresponding to $\nabla^\pm =\nabla\pm\psi$ and ${}^\tau{}\nabla$ corresponding to a connection $\nabla$ was defined in Eq.~\eqref{tau-conn}.
 \end{deff}

\begin{rem} \label{nablai}
Since the difference of any two connections on a bundle $A$ is a section in $T^*M\otimes \End(A)$, one may consider $\nabla^+=:\nabla^1$ and $\nabla^-=:\nabla^2$ as two independent connections, with $\nabla^i$ acting on the $i$-th slot or factor of $\Phi \in \Gamma(T^*M \otimes T^*M)$, where here $i\in \{1,2\}$. So, the data $(A,\nabla,\psi)$ can be also replaced by $(A,\nabla^1,\nabla^2)$ in the above definition. The reason for using the parametrization with $\psi$ becomes clear from the following:
\end{rem}

\begin{rem} \label{symskew}
Let us decompose $\Phi$ into the sum of symmetric and skew-symmetric parts, $\Phi=\Phi^{\mathrm{sym}}+\Phi^{\mathrm{skew}}$.
Using the defining equations, it is easy to verify that for all $X,Y \in \Gamma(TM)$ and $s \in \Gamma(A)$ one has
 $$ \left({}^\tau{}\nabla^{comb}_s \Phi \right)(X,Y) = \left({}^\tau{}\nabla_s \Phi \right)(X,Y) - \Phi( \langle \rho(\psi s) , X \rangle , Y) +
\Phi(X,\langle \rho(\psi s), Y\rangle) \: ;$$
here $\psi s \equiv \psi(s) \in \Omega^1 (M,A)$ and
$\rho (\psi s) \in \Gamma(T^*M \otimes TM) \cong \Gamma(\End (TM))$. Note that ${}^\tau{}\nabla^{comb}$ does not map symmetric 2-tensors into symmetric ones, while ${}^\tau{}\nabla$ does. As a consequence, Eq.~\eqref{Bilinear-Lie-Killing} mixes $\Phi^{\mathrm{sym}}$ and $\Phi^{\mathrm{skew}}$; in detail, it decomposes into
\bqa\label{Bilinear_symm}
{}^\tau{}\nabla_s ({\Phi^{\mathrm{sym}}}) &=&\mathrm{Sym}\langle \rho(\psi s), \Phi^{\mathrm{skew}}\rangle   \: , \\ \nonumber
{}^\tau{}\nabla_s (\Phi^{\mathrm{skew}}) &=& \mathrm{Alt}\langle \rho(\psi s), \Phi^{\mathrm{sym}}\rangle
\eqa
 for all $s\in\Gamma (A)$. Here the contraction is taken by the natural pairings $TM\otimes \Lambda^2(T^*M)\to T^*M$
and $TM\otimes \mathrm{Sym}^2(T^*M)\to T^*M$, respectively, extending the pairing between $TM$ and $T^*M$ by the Leibniz rule.\footnote{If, alternatively, the contraction of $TM$ is defined to be taken with the first factor of the respective 2-tensor only, there is a factor of 2 to insert on the r.h.s.~of each equation.} Using Eq.~\eqref{tau-conn}, this in turn can be rewritten in the form of Eqs.~\eqref{Lie3} with the upper sign ($\Phi^{\mathrm{sym}} := g$,  $\Phi^{\mathrm{skew}}:=B$) as found in the string-physics applications (cf.~also \cite{withyoungpeople}). More explicitly, let $(e_a)_{a=1}^{\mathrm{rk}A}$ be a local basis of sections, such that
$\nabla e_a=\omega_a^b e_b$, $\psi e_a=\psi_a^b e_b$, and $\rho_a=\rho(e_a)$, then the Eqs.~(\ref{Bilinear_symm}) read as follows:
\beq\label{Generalized_symmetric}
\cL_{\rho_a}g = \omega_a^b \vee\iota_{\rho_b} g + \psi_a^b \vee \iota_{\rho_b} B \, ,\\ \nonumber 
\cL_{\rho_a}B = \omega_a^b \wedge\iota_{\rho_b} B + \psi_a^b \wedge \iota_{\rho_b} g \, .
\eeq
We will come back to this example again.
\end{rem}

\noindent The following construction extends that of Example \ref{reg_foliation} (cf.~\cite{Kotov-Strobl14}). Before stating the proposition to be proven, we set up some noations and terminology. A bilinear form $b$ on a vector space $T$
uniquely corresponds to a linear map $\hat b\colon T\to T^*$ by $b(v_1, v_2) =\langle \hat b (v_1), v_2\rangle$
for all $v_1, v_2\in V$; the bilinear form is non-degenerate if and only if $\hat b$ is invertible.
It is obvious that every bilinear form with a positive (or negative) definite symmetric part is non-degenerate.
Whenever the inverse map $\hat b^{-1}$ exists,  
it determines a bilinear form on $T^*$ by means of
$b^{-1}(\lambda_1, \lambda_2):=\langle \lambda_1, \hat b^{-1} (\lambda_2)\rangle$ for $\lambda_1, \lambda_2\in T^*$. We call a map $\lambda \colon (T_1, b_1) \to (T_2,b_2)$ such that $b_1(u,v) = b_2(\lambda(u),\lambda(v))$ for all $u,v \in T_1$ a \emph{generalized isometry}. With the above definitions and for a non-degenerate $b$,  $\hat b \colon (T,b) \to (T^*,b^{-1})$ is such a generalized isometry. For every vector subspace $V\subset T$ the \emph{left-orthogonal} subspace is defined as
$V^\perp =\{ v\in T \mid b(v,V)=0\}$. It follows immediately that $V$ and $V^\perp$ are complimentary if and only if the restriction of $b$ onto $V$ is non-degenerate; e.g.~this is always satisfied if the symmetric part of $b$ is positive (or negative) definite. An easy computation shows that an invertible $\hat b$ gives rise to a generalized isometry between
$(V^\perp, b|_{V^\perp})$ and $(\mathrm{Ann}(V), b^{-1}|_{\mathrm{Ann}(V)})$, where
$\mathrm{Ann}(V)\subset T^*$ is the annihilator of $V$.

\begin{proposition}\label{Bilinear_quotient} Given a (regular) foliation ${\cF}$ of a
smooth manifold $M$ equipped with a non-degenerate bilinear form $\Phi$, then its co-normal bundle $N^*\cF$
with the canonically induced bilinear form is invariant with respect to leaf-preserving diffeomorphisms on $M$ if and only if
there exist connections $\nabla^\pm$ on $A=T\cF$ such that (\ref{Bilinear-Lie-Killing}) holds true.
\end{proposition}

\begin{cor}
If the foliated manifold has, in addition, a smooth quotient $Q=M/\cF$,
then $T^*Q$ can be equipped with a bilinear form such that for any point $p \in M$ the natural linear isomorphism $N^*_p M \to T_{[p]}^*Q$
is a generalized isometry.
\end{cor}

\noindent{\bf Proof} {(of Proposition \ref{Bilinear_quotient}){\bf .}} Let $\beta_1$, $\beta_2$ be sections of $N^*\cF$, the annihilator of $T\cF$, which are transversally invariant,
i.e.~invariant under the action of ``vertical'' vector fields (vector fields parallel to $T\cF$; this allows to identify $\beta_1$ and $\beta_2$ with 1-forms on the quotient space $Q$ whenever it is smooth).
As soon as we prove that $\Phi^{-1}(\beta_1, \beta_2)$ is transversally invariant if and and only if (\ref{Bilinear-Lie-Killing}) is fulfilled for some $\nabla^\pm$, the assertion of Proposition \ref{Bilinear_quotient} follows.
 Indeed, since the above 1-forms are invariant under the Lie derivative along any section $s$ of $A=T\cF$,
 we must require that $\cL_s \left(\Phi^{-1}\right)(\beta_1, \beta_2)=0$. Therefore the restriction of $\cL_s \left(\Phi^{-1}\right)$ to $N^*\cF$ vanishes at every point, which is true if and only if $\cL_s \left(\Phi^{-1}\right)\in \Gamma (T\cF\otimes TM + TM\otimes T\cF)$.
Taking into account that $\Phi$ is non-degenerate, we get that for any locally defined basis of vector fields
 $(\rho_a)_{a=1}^{r}$ of $T{\cF}$, $r=\dim \cF$, there exist two local $r\times{r}$ matrices $\omega^\pm$ such that the following identity holds true:
\begin{equation} \label{decomp}
{\cL}_{\rho_a} \left(\Phi^{-1}\right) +  \rho_b\otimes \langle \left(\omega^+\right)_a^b, \Phi^{-1}\rangle  
+ \langle\Phi^{-1}, \left(\omega^-\right)_a^b\rangle \otimes \rho_b =0\,,
\end{equation}
where the second and the third terms of the l.h.s.~are left and right contractions of a contravariant $2$-tensor with a 1-form, respectively.
By formula (\ref{tau-conn}), extended to sections of $TM\otimes TM$, we conclude that
the identity (\ref{decomp}) is equivalent to the existence
of local connections on $A$ such that ${}^\tau{}\nabla^{comb} (\Phi^{-1})=0$ and thus (\ref{Bilinear-Lie-Killing}) is satisfied.
Using a partition of unity, we obtain global connections on $A$ with the same property.
$\square$

\begin{rem}\label{remark_invertible_on_the_quotient}
If the restriction of $\Phi$ onto $T\cF$ is non-degenerate, then so is the induced bilinear form on $T\cF^\perp$ and, by the canonical generalized isometry, on
$N^*\cF$. Hence the bilinear form on $T^*Q$, obtained under the assumptions of Proposition \ref{Bilinear_quotient}, has an inverse, which we will denote by $\Phi_Q$.
\end{rem}

\begin{example}\label{Generalized-Lie-Killing_example} A \emph{generalized Riemannian structure} on an $n-$dimensional manifold $M$ is a rank $n$ subbundle of the exact Courant algebroid
$TM\oplus T^*M$ on which the inner product is positive definite; this construction has been used in several occasions already, but explicitly introduced in particular by N.~Hitchin (\cite{Hitchin03}, cf. also \cite{Kotov-Strobl10} for an exposition of generalized geometry). Generalized Riemannian structures are in one-to-one correspondence with bilinear forms
 $\Phi=g+B$, where $g$ is a Riemannian metric tensor and $B$ is a skew-symmetric 2-form; the correspondence is given by the graph of $\Phi$ considered as a bundle map $TM\to T^*M$. Then $(A,\nabla^\pm)$ is compatible with $(M,g,B)$ if
${}^\tau{}\nabla^{comb} (g+B)=0$, rewritten in other ways also in Remark \ref{symskew} above.
Let us notice that the restriction of a generalized Riemannian structure on any subbundle of $TM$ is non-degenerate
as its symmetric part is positive definite. Thus under the assumptions of Proposition \ref{Bilinear_quotient}, Remark \ref{remark_invertible_on_the_quotient} permits us to conclude that we obtain
a generalized Riemannian structure $\Phi_Q \equiv g_Q + B_Q$ on the quotient space $Q=M/{\cF}$.
\end{example}

If a bilinear form is totally skew-symmetric, $\Phi = B \in \Omega^2(M)$, and setting $\psi=0$, the equations (\ref{Generalized_symmetric})
reduce to ${}^\tau{}\nabla (B)=0$, a compatibility condition similar in spirit to (\ref{Lie-Killing1}); locally this equation takes the form
$\cL_{\rho_a} B = \omega_a^b \wedge \iota_{\rho_a} B$ (cf.~Remark \ref{symskew}). We now consider $B$ to be a symplectic form:

\begin{deff}\label{symplectic-Lie_def} Let $(A,\rho,[\cdot,\cdot])$  be a Lie algebroid over a symplectic manifold $(M, \Omega)$ and
$\nabla$ be a connection on $A$.
Then $(A,\nabla)$ and $(M,\Omega)$ are called \emph{compatible}, if
\begin{equation}\label{symplectic-Lie}
{}^\tau{}\nabla (\Omega)=0\,.
\end{equation}
 \end{deff}

\begin{example}\label{symplectic_quotient}
Given a (regular) foliation $\cF$ on a symplectic manifold $(M,\Omega)$ with a smooth quotient $Q=M/\cF$ and a connection $\nabla$ on $A=T\cF$,
such that the compatibility condition \eqref{symplectic-Lie} is fulfilled, we immediately get a canonical bivector field on the quotient space.
This follows by the same method as in Proposition \ref{Bilinear_quotient}.
The bivector field on $Q$ obtained above is clearly Poisson and so is the quotient map.
In addition, the quotient Poisson structure is symplectic if and only if the restriction of $\Omega$ on the fibers of $\cF$ is non-degenerate.
\end{example}

\noindent Defintion \ref{symplectic-Lie_def} and
Example \ref{symplectic_quotient} admit a straightforward generalization to the Poisson case:  a Lie algebroid $A$ with a connection $\nabla$ over
$M$ is compatible with a Poisson structure $\cP$ if ${}^\tau{}\nabla (\cP)=0$. Under the assumptions of Example \ref{symplectic_quotient}, we get a canonical
Poisson structure on the quotient space $Q=M/\cF$ such that the quotient map is Poisson.

\begin{rem}
For a regular foliation $\cF$, transversal invariance of a tensor field is \emph{equivalent} to the existence of a connection $\nabla$ on $T\cF$ such that ${}^\tau{}\nabla$ annihilates this tensor field. However, for non-regular foliations, annihilation by some ${}^\tau{}\nabla$ is a \emph{stronger} requirement than transversal invariance: The standard metric $g$ on $\R^2$ is transversally invariant with respect to the singular foliation determined by the Lie algebroid from Proposition \ref{rho0}. However, we proved in Proposition \ref{rho0} that there is no connection $\nabla$ such that Eq.~\eqref{Lie-Killing1} is satisfied.
\end{rem}

\section*{Appendix A: Flat Killing Lie algebroids, simple examples and facts}

\noindent A typical example of a flat Killing Cartan Lie algebroid, i.e.~a Killing Cartan Lie algebroid as defined in Definition \ref{CartanRiemann-def} such that $\nabla$ is flat but which is not just an action Lie algebroid mapping into isometries by means of $\rho$, but locally so only (cf.~also Example \ref{example_flat_CartanKilling}), is the following one:\footnote{This example was suggested to us by A.~Weinstein.}
\begin{example} Consider the unit square  $[0,1] \times [0,1] \in \R^2$ with an identification of opposite sides so as to yield a Klein bottle $M$. Equipping $M$ with the natural metric $g$ and its
Levi-Civita connection $\nabla$, we get another example of a flat Killing Lie algebroid. Similarly, the foliation of $\R^2$ by vertical lines induces a foliation $\cF$ on $M$ and ($T\cF$,$\nabla$,$g$) forms a sub-Killing Lie algebroid of the previous one. In both cases, locally constant sections correspond to local Killing vectors, while some of them do not extend to globally constant (non-zero) sections due to the non-triviality of the bundles.  \label{Klein-bottle}
\end{example}
\noindent For the same reason, $TM$ and $T\cF$ cannot be action Lie algebroids in the last example, while locally they are, corresponding to the fact that some of the local isometries of the Klein bottle do not extend to global ones.

\vskip 3mm \noindent Note that one cannot expect to necessarily always find \emph{all} Killing vectors as coming from constant sections, also not locally. Here two examples, both of which are in fact also (flat Killing) Cartan Lie algebroids:
\begin{example}\label{example_Mink}
Let $M=\R^n$ with its standard flat metric $g$, $A$ the standard Lie algebroid $A := TM \cong \R^n \times \R^n$, equipped with its standard flat connection of a vector space. It shows that there are Killing vectors which do not arise from flat sections: While the generators of translations, $\partial_{i}$, are constantly covariant and Killing vectors by the above argument, the generators of rotations, $m_{ij} := x^i \partial_j - x^j \partial_i$, $i \neq j$, are Killing vectors of $g$
 without being covariantly constant. In fact, there is \emph{no} (globally defined) connection $\nabla$ on  $TM$ making these generators covariantly constant.
\end{example}
\begin{example} \label{example_sphere}
Let $M_0 = S^2 \subset \R^3$ equipped with its standard metric of constant curvature. Denote by $e_1$, $e_2$, and $e_3$ a basis of the three-dimensional isometry Lie algebra $so(3)$. Consider a region $M \subset M_0$ where $e_1$ and $e_2$ are non-zero everywhere. Let $A=TM$ be the standard Lie algebroid over this $M$. Define a connection on $A$ by requiring $\nabla e_1=0=\nabla e_2$. Now since on $M$ one has $e_3 = f_1 e_1 + f_2 e_2$ for non-constant functions $f_1$ and $f_2$, $\nabla e_3 \neq 0$, while still $\rho(e_3)=[\rho(e_1),\rho(e_2)]$ is a Killing vector field.
\end{example}

\noindent In both of these examples, the bundle has smaller rank than the dimension of isometry group of the Riemannian base manifold. But the situation in this context here can be also reversed easily: There is a simple procedure to construct flat Killing Lie algebroids of higher rank from any given: Let $(A_0,\nabla_0,g)$ be a flat Killing Lie algebroid over a manifold $M$. Choose any Lie algebroid extension $A$,
$$ 0 \to B \to A \to A_0 \to 0  \, , $$
of $A_0$ by a bundle of Lie algebras $B$ over $M$. Assume that $A$ admits a flat connection $\nabla$  projecting to $\nabla_0$ on $A_0$. Then the triple $(A,\nabla,g)$ also defines a flat Killing Lie algebroid. Note that here the image of $\rho$ and of  $\rho_0$ coincide at each $x \in M$.

\vskip 3mm \noindent
In general, the relation of a flat Killing Lie algebroid to an action Lie algebroid with identically induced foliation can be intricate. What always holds true, however, is stated in the following
\begin{proposition}\label{Lemma2}
Let $(A,\nabla,g)$ be a flat Killing Lie algebroid over $M$, i.e.~a Killing Lie algebroid where $\nabla$ is flat. Then every point $x \in M$ permits a neighborhood $x \ni U \subset M$ over which there exists a finite-rank action Lie algebroid the induced singular foliation of which is identical to the one of $A|_U$. This action Lie algebroid is canonically a flat Killing Cartan Lie algebroid. \label{lemmaflat}
\end{proposition}
\noindent \proof For every $x \in M$ we may choose a neighborhood with a covariantly constant basis $(e_a)_{a=1}^{\rk A}$ of sections $e_a \in \Gamma(A|_U)$ satisfying $\nabla e_a =0$. Its image by $\rho$ provides a set of Killing vectors which,  over $C^\infty(U)$, generates the image of $\rho|_U$---the integral surfaces of this (possibly singular) distribution provides the (possibly singular) foliation of $U$. Taken over $\R$, the same set of Killing vectors generates a finite-dimensional Lie algebra $\g_U$ acting on $U$. $U \times\g_U$ is the searched-for action Lie algebroid.  $\square$

\begin{rem}
In both examples Example \ref{example_Mink} and Example \ref{example_sphere}, the action Lie algebroid constructed along the above lines glues together to globally flat Killing Cartan Lie algebroids $M \times \mathrm{iso(g)}$ where $\mathrm{iso(g)}$ is the full isometry Lie algebra of $g$. In Example  \ref{example_sphere}, e.g., $M \times \mathrm{so}(3)$ (equipped with its canonical flat connection).
There now is a Lie algebroid morphism from this action Lie algebroid into the original (flat Killing Cartan) Lie algebroid $A=TM$, by mapping the standard $so(3)$-basis to $e_1$, $e_2$, $e_3$. However, this morphism does not respect the connection $\nabla$, and thus is not a morphism in the category of Killing Lie algebroids. In Example \ref{Klein-bottle}, on the other hand, we found a local isomorphism as (flat Cartan) Killing Lie algebroids.
\end{rem}


\end{document}